\documentclass[11pt]{article}
\usepackage[]{amsmath,amssymb}
\usepackage{cite}
\newtheorem{theorem}{Theorem}[section]
\newtheorem{lemma}[theorem]{Lemma}
\newtheorem{proposition}[theorem]{Proposition}
\newtheorem{corollary}[theorem]{Corollary}
\newtheorem{exAux}[theorem]{Example}

\newtheorem{Def}[theorem]{Definition}
\newenvironment{definition}{\begin{Def} \rm}{\end{Def}}
\newtheorem{Note}[theorem]{Note}
\newenvironment{note}{\begin{Note} \rm}{\end{Note}}
\newtheorem{Problem}[theorem]{Problem}
\newenvironment{problem}{\begin{Problem} \rm}{\end{Problem}}
\newtheorem{Rem}[theorem]{Remark}

\newtheorem{Not}[theorem]{Notation}

\newtheorem{Conj}[theorem]{Conjecture}
\newenvironment{conjecture}{\begin{Conj}}{\end{Conj}}
\newtheorem{Ass}[theorem]{Assumption}

\newenvironment{proof}{\medskip\noindent{\bf Proof.\ }}{\qed\medskip}
\newenvironment{proofof}[1]{\medskip\noindent{\bf Proof  of {#1}.\ 
}}{\qed\medskip}
\newcommand{\qed}{\hfill\mbox{$\Box$\qquad\qquad}}

\newcommand{\Mat}[1]{\text{\rm Mat}_{#1}(\mathbb{K})}

\renewcommand{\b}[1]{\langle #1 \rangle}
\newcommand{\tr}{\text{\rm tr}}
\renewcommand{\d}{\downarrow}
\newcommand{\D}{\Downarrow}

\renewcommand{\th}{\theta}

\renewcommand{\indent}{\hspace{6mm}}

\begin{document}
\thispagestyle{empty}

\begin{center}
\LARGE \bf
\noindent
Sharp tridiagonal pairs
\end{center}

\smallskip

\begin{center}
\Large
Kazumasa Nomura and 
Paul Terwilliger\footnote{This author gratefully acknowledges 
support from the FY2007 JSPS Invitation Fellowship Program
for Reseach in Japan (Long-Term), grant L-07512.}
\end{center}

\smallskip

\begin{quote}
\small 
\begin{center}
\bf Abstract
\end{center}

\indent
Let $\mathbb{K}$ denote a field and let $V$ denote a vector space
over $\mathbb{K}$ with finite positive dimension.
We consider a pair of $\mathbb{K}$-linear transformations $A:V \to V$
and $A^*:V \to V$ that satisfies the following conditions:
(i)
each of $A,A^*$ is diagonalizable;
(ii)
there exists an ordering $\{V_i\}_{i=0}^d$ of the eigenspaces of 
$A$ such that
$A^* V_i \subseteq V_{i-1} + V_{i} + V_{i+1}$ for $0 \leq i \leq d$,
where $V_{-1}=0$ and $V_{d+1}=0$;
(iii)
there exists an ordering $\{V^*_i\}_{i=0}^\delta$ 
of the eigenspaces of $A^*$ such that
$A V^*_i \subseteq V^*_{i-1} + V^*_{i} + V^*_{i+1}$ for $0 \leq i \leq \delta$,
where $V^*_{-1}=0$ and $V^*_{\delta+1}=0$;
(iv) 
there is no subspace $W$ of $V$ such that
$AW \subseteq W$, $A^* W \subseteq W$, $W \neq 0$, $W \neq V$.
We call such a pair a {\em tridiagonal pair} on $V$.
It is known that $d=\delta$ and for $0 \leq i \leq d$
the dimensions of $V_i$, $V_{d-i}$, $V^*_i$, $V^*_{d-i}$ coincide.
We say the pair $A,A^*$ is {\em sharp} whenever $\dim V_0=1$.
A conjecture of Tatsuro Ito and the second author
states that if $\mathbb{K}$ is algebraically closed then $A,A^*$ is sharp.
In order to better understand and eventually prove the conjecture,
in this paper we begin a systematic study of the sharp tridiagonal pairs.
Our results are summarized as follows. 
Assuming $A,A^*$ is sharp and using the data 
$\Phi=(A; \{V_i\}_{i=0}^d; A^*; \{V^*_i\}_{i=0}^d)$
we define a finite sequence of scalars called the parameter array. 
We display some equations that show the geometric
significance of the parameter array. 
We show how the parameter array is affected if $\Phi$ is
replaced by
$(A^*; \{V^*_i\}_{i=0}^d; A; \{V_i\}_{i=0}^d)$
or
$(A; \{V_{d-i}\}_{i=0}^d; A^*; \{V^*_i\}_{i=0}^d)$
or
$(A; \{V_{i}\}_{i=0}^d; A^*; \{V^*_{d-i}\}_{i=0}^d)$.
We prove that if the isomorphism class of $\Phi$ is determined by
the parameter array then there exists a nondegenerate symmetric
bilinear form $\b{\;,\,}$ on $V$ such that
$\b{Au,v}=\b{u,Av}$ and $\b{A^*u,v}=\b{u,A^*v}$ for all $u,v \in V$.
\end{quote}

\section{Tridiagonal pairs}

\indent
Throughout the paper $\mathbb{K}$  denotes a field and
$V$ denotes a vector space over $\mathbb{K}$ with finite
positive dimension.

\medskip

We begin by recalling the notion of a tridiagonal pair.
We will use the following terms.
Let $\text{End}(V)$ denote the $\mathbb{K}$-algebra 
consisting of all $\mathbb{K}$-linear transformations from $V$ to $V$.
For $A \in \text{End}(V)$ and for a subspace $W \subseteq V$,
we call $W$ an {\em eigenspace} of $A$ 
whenever $W \neq 0$ and there exists $\theta \in \mathbb{K}$
such that $W=\{v \in V \,|\, Av=\theta v\}$; 
in this case $\th$ is
the {\em eigenvalue} of $A$ associated with $W$.
We say $A$ is {\em diagonalizable} whenever $V$ is spanned by
the eigenspaces of $A$.
We now recall the notion of a tridiagonal pair.

\medskip

\begin{definition} \cite{ITT}     \label{def:TDpair}   \samepage
By a {\em tridiagonal pair} on $V$ we mean an ordered pair
of elements $A,A^*$ taken from $\text{End}(V)$ that satisfy (i)--(iv) below:
\begin{itemize}
\item[(i)] 
Each of $A,A^*$ is diagonalizable.
\item[(ii)] 
There exists an ordering $\{V_i\}_{i=0}^d$ of the eigenspaces of 
$A$ such that
\begin{equation}              \label{eq:Astrid}
  A^* V_i \subseteq V_{i-1} + V_{i} + V_{i+1}  \qquad\qquad (0 \leq i \leq d),
\end{equation}
 where $V_{-1}=0$ and $V_{d+1}=0$.
\item[(iii)] 
There exists an ordering $\{V^*_i\}_{i=0}^\delta$ 
of the eigenspaces of $A^*$ such that
\begin{equation}           \label{eq:Atrid}
A V^*_i \subseteq V^*_{i-1} + V^*_{i} + V^*_{i+1} \qquad\qquad
            (0 \leq i \leq \delta),
\end{equation}
where $V^*_{-1}=0$ and $V^*_{\delta+1}=0$.
\item[(iv)] 
There is no subspace $W$ of $V$ such that
$AW \subseteq W$, $A^* W \subseteq W$, $W \neq 0$, $W \neq V$.
\end{itemize}
We say the pair $A,A^*$ is {\em over} $\mathbb{K}$.
\end{definition}

\begin{note}      \label{note:star}        \samepage
It is a common notational convention to use $A^*$ to represent the
conjugate-transpose of $A$. We are not using this convention.
In a tridiagonal pair $A,A^*$ the linear transformations $A$ and
$A^*$ are arbitrary subject to (i)--(iv) above.
\end{note}

\medskip
We refer the reader to
\cite{AC,AC2,AC3,Bas,ITT,IT:shape,IT:uqsl2hat,IT:Krawt,N:aw,N:refine,N:height1,
NT:tde,Vidar} for background on tridiagonal pairs.
See 
\cite{BI,BT:Borel,BT:loop,Bow,Ca,CaMT,CaW,Egge,F:RL,H:tetra,HT:tetra,
IT:non-nilpotent,IT:tetra,
IT:inverting,IT:drg,IT:loop,ITW:equitable,Leo,Mik,Mik2,
P,R:multi,R:6j,Suz,T:sub1,T:sub3,
T:qSerre,T:Kac-Moody,Z}
for related topics.

\medskip

Let $A,A^*$ denote a tridiagonal pair on $V$, 
as in Definition \ref{def:TDpair}.
By \cite[Lemma 4.5]{ITT} the integers $d$ and $\delta$ from (ii), (iii) are equal;
we call this common value the {\em diameter} of the pair. 
By \cite[Corollary 5.7]{ITT}, for $0 \leq i \leq d$ the spaces $V_i$, $V^*_i$
have the same dimension; we denote
this common dimension by $\rho_i$. 
By the construction $\rho_i \neq 0$.
By \cite[Corollaries 5.7, 6.6]{ITT}
the sequence $\{\rho_i\}_{i=0}^d$ is symmetric and unimodal;
that is $\rho_i=\rho_{d-i}$ for $0 \leq i \leq d$ and
$\rho_{i-1} \leq \rho_i$ for $1 \leq i \leq d/2$.
We call the sequence $\{\rho_i\}_{i=0}^d$ the {\em shape}
of $A,A^*$.
The following special case has received a lot of attention. 
By a {\em Leonard pair} we mean a tridiagonal pair with shape
$(1,1,\ldots,1)$ \cite[Definition 1.1]{T:Leonard}.
There is a natural correspondence between the Leonard pairs and a
family of orthogonal polynomials consisting of the $q$-Racah polynomials
and their relatives \cite{T:qRacah}.
This family coincides with the terminating branch of the Askey scheme
\cite{Koe}.
See 
\cite{Cur:mlt,Cur:spinLP,H,M:LT,NT:balanced,NT:formula,NT:det,NT:mu,
NT:span,NT:switch,NT:affine,NT:maps,T:Leonard,T:24points,T:canform,T:intro,
T:intro2,T:split,T:array,T:qRacah,T:survey,TV,V,V:AW}
for more information about Leonard pairs.
Our point of departure in this paper is the following conjecture.

\medskip

\begin{conjecture}  {\rm \cite{IT:Krawt}}  \label{conj:closed}   \samepage
Let $A,A^*$ denote a tridiagonal pair over an algebraically
closed field. Then $\rho_0=1$ where $\{\rho_i\}_{i=0}^d$ is the shape of $A,A^*$.
\end{conjecture}

\begin{note}
Conjecture \ref{conj:closed} has been proven for some special cases
that are described as follows. Referring to Conjecture \ref{conj:closed},
there is a parameter $q$ associated with $A,A^*$ 
that is used to describe the eigenvalues; we discuss
$q$ above Lemma \ref{lem:formula}. In \cite{IT:aug} Ito and the second author
prove Conjecture \ref{conj:closed} assuming $q$ is not a root of unity.
There is a family of tridiagonal pairs with $q=1$
that are said to have {\em Krawtchouk type} \cite{IT:Krawt}. 
In \cite{IT:Krawt} Ito and the second author prove 
Conjecture \ref{conj:closed} assuming $A,A^*$
has Krawtchouk type.
\end{note}

\medskip

In order to better understand and eventually prove Conjecture \ref{conj:closed},
in this paper we begin a systematic study of those tridiagonal pairs that
satisfy its conclusion. We start with a definition.

\medskip

\begin{definition}       \samepage
A tridiagonal pair $A,A^*$ is said to be {\em sharp} whenever
$\rho_0=1$, where $\{\rho_i\}_{i=0}^d$ is the shape of $A,A^*$.
\end{definition}

\section{Tridiagonal systems}

\indent
When working with a tridiagonal pair, it is often convenient to consider
a closely related object called a tridiagonal system.
To define a tridiagonal system, we recall a few concepts from linear
algebra.
Let $A$ denote a diagonalizable element of $\text{End}(V)$.
Let $\{V_i\}_{i=0}^d$ denote an ordering of the eigenspaces of $A$
and let $\{\th_i\}_{i=0}^d$ denote the corresponding ordering of
the eigenvalues of $A$.
For $0 \leq i \leq d$ let $E_i:V \to V$ denote the linear transformation
such that $(E_i-I)V_i=0$ and $E_iV_j=0$ for $j \neq i$ $(0 \leq j \leq d)$.
Here $I$ denotes the identity of $\text{End}(V)$.
We call $E_i$ the {\em primitive idempotent} of $A$ corresponding to $V_i$
(or $\th_i$).
Observe that
(i) $\sum_{i=0}^d E_i = I$;
(ii) $E_iE_j=\delta_{i,j}E_i$ $(0 \leq i,j \leq d)$;
(iii) $V_i=E_iV$ $(0 \leq i \leq d)$;
(iv) $A=\sum_{i=0}^d \th_iE_i$.
Moreover
\begin{equation}         \label{eq:defEi}
  E_i=\prod_{\stackrel{0 \leq j \leq d}{j \neq i}}
          \frac{A-\th_jI}{\th_i-\th_j}.
\end{equation}
We note that $\{E_i\}_{i=0}^d$ is a basis for the $\mathbb{K}$-subalgebra
of $\text{End}(V)$ generated by $A$.

Now let $A,A^*$ denote a tridiagonal pair on $V$.
An ordering of the eigenspaces of $A$ (resp. $A^*$)
is said to be {\em standard} whenever it satisfies 
\eqref{eq:Astrid} (resp. \eqref{eq:Atrid}). 
We comment on the uniqueness of the standard ordering.
Let $\{V_i\}_{i=0}^d$ denote a standard ordering of the eigenspaces of $A$.
Then the ordering $\{V_{d-i}\}_{i=0}^d$ is standard and no other ordering
is standard.
A similar result holds for the eigenspaces of $A^*$.
An ordering of the primitive idempotents of $A$ (resp. $A^*$)
is said to be {\em standard} whenever
the corresponding ordering of the eigenspaces of $A$ (resp. $A^*$)
is standard.
We now define a tridiagonal system.

\medskip

\begin{definition}          \label{def:TDsystem}  \samepage
By a {\em tridiagonal system} on $V$ we mean a sequence
\[
 \Phi=(A;\{E_i\}_{i=0}^d;A^*;\{E^*_i\}_{i=0}^d)
\]
that satisfies (i)--(iii) below.
\begin{itemize}
\item[(i)]
$A,A^*$ is a tridiagonal pair on $V$.
\item[(ii)]
$\{E_i\}_{i=0}^d$ is a standard ordering
of the primitive idempotents of $A$.
\item[(iii)]
$\{E^*_i\}_{i=0}^d$ is a standard ordering
of the primitive idempotents of $A^*$.
\end{itemize}
We say $\Phi$ is {\em over} $\mathbb{K}$.
\end{definition}

\medskip

The following result is immediate from lines \eqref{eq:Astrid}, \eqref{eq:Atrid}
and Definition \ref{def:TDsystem}.

\medskip

\begin{lemma}    \label{lem:trid}     \samepage
Let $(A;\{E_i\}_{i=0}^d;A^*;\{E^*_i\}_{i=0}^d)$ denote a tridiagonal system.
Then for $0 \leq i,j \leq d$ the following {\rm (i)}, {\rm (ii)} hold.
\begin{itemize}
\item[\rm (i)]
$E_i A^* E_j = 0\;\;$ if  $|i-j|>1$.
\item[\rm (ii)]
$E^*_i A E^*_j =0\;\;$ if $|i-j|>1$.
\end{itemize}
\end{lemma}

\begin{definition}        \label{def}  \samepage
Let $\Phi=(A;\{E_i\}_{i=0}^d;A^*$; $\{E^*_i\}_{i=0}^d)$ 
denote a tridiagonal system on $V$.
For $0 \leq i \leq d$ let $\th_i$ (resp. $\th^*_i$)
denote the eigenvalue of $A$ (resp. $A^*$)
associated with the eigenspace $E_iV$ (resp. $E^*_iV$).
We call $\{\th_i\}_{i=0}^d$ (resp. $\{\th^*_i\}_{i=0}^d$)
the {\em eigenvalue sequence}
(resp. {\em dual eigenvalue sequence}) of $\Phi$.
We observe $\{\th_i\}_{i=0}^d$ (resp. $\{\th^*_i\}_{i=0}^d$) are
mutually distinct and contained in $\mathbb{K}$.
By the {\em shape} of $\Phi$ we mean the shape
of the tridiagonal pair $A,A^*$.
We say $\Phi$ is {\em sharp} whenever the tridiagonal pair $A,A^*$ is
sharp.
\end{definition}

\medskip

We have a comment.

\medskip

\begin{lemma} {\rm \cite[Theorem 11.1]{ITT}}  \label{lem:indep} \samepage
Let $\Phi$ denote a tridiagonal system 
with eigenvalue sequence $\{\th_i\}_{i=0}^d$ and 
dual eigenvalue sequence $\{\th^*_i\}_{i=0}^d$.
Then the expressions
\begin{equation}         \label{eq:beta}
 \frac{\th_{i-2}-\th_{i+1}}{\th_{i-1}-\th_i},
\qquad\qquad
 \frac{\th^*_{i-2}-\th^*_{i+1}}{\th^*_{i-1}-\th^*_i}
\end{equation}
are equal and independent of $i$ for $2 \leq i \leq d-1$.
\end{lemma}

\medskip

When discussing tridiagonal systems we will use the following
notational convention.
Let $\lambda$ denote an indeterminate and let $\mathbb{K}[\lambda]$
denote the $\mathbb{K}$-algebra consisting of all polynomials
in $\lambda$ that have coefficients in $\mathbb{K}$.
With reference to Definition \ref{def}
for $0 \leq i \leq d$ we define the following polynomials in
$\mathbb{K}[\lambda]$:
\begin{align*}
 \tau_i &= 
  (\lambda-\th_0)(\lambda-\th_1)\cdots(\lambda -\th_{i-1}), \\
 \eta_i &=
  (\lambda-\th_d)(\lambda-\th_{d-1})\cdots(\lambda-\th_{d-i+1}),  \\
 \tau^*_i &=
  (\lambda-\th^*_0)(\lambda-\th^*_1)\cdots(\lambda-\th^*_{i-1}), \\
 \eta^*_i &=
  (\lambda-\th^*_d)(\lambda-\th^*_{d-1})\cdots(\lambda-\th^*_{d-i+1}).
\end{align*}
Note that each of $\tau_i$, $\eta_i$, $\tau^*_i$, $\eta^*_i$ is
monic with degree $i$.
By \eqref{eq:defEi}, for $0 \leq i \leq d$ we have
\begin{align*}
 E_i &=\frac{\tau_i(A)\eta_{d-i}(A)}
            {\tau_i(\th_i)\eta_{d-i}(\th_i)}, 
 &
 E^*_i &=\frac{\tau^*_i(A^*)\eta^*_{d-i}(A^*)}
              {\tau^*_i(\th^*_i)\eta^*_{d-i}(\th^*_i)}.
\end{align*}
In particular
\begin{align}
  E_0 &= \frac{\eta_d(A)}{\eta_d(\th_0)}, &
  E_d &= \frac{\tau_d(A)}{\tau_d(\th_d)},      \label{eq:E0}   \\
  E^*_0 &= \frac{\eta^*_d(A^*)}{\eta^*_d(\th^*_0)}, &
  E^*_d &= \frac{\tau^*_d(A^*)}{\tau^*_d(\th^*_d)}.      \label{eq:Es0}
\end{align}

\medskip

We mention a result for future use.

\medskip

\begin{lemma}{\rm \cite[Proposition 5.5]{NT:mu}}     \samepage
With reference to Definition {\rm \ref{def}},
\begin{equation}         \label{eq:etad}
 \eta_d=\sum_{i=0}^d \eta_{d-i}(\th_0)\tau_i,
\qquad\qquad
 \eta^*_d=\sum_{i=0}^d \eta^*_{d-i}(\th^*_0)\tau^*_i.
\end{equation}
\end{lemma}

\section{Isomorphisms of tridiagonal systems}

\indent
Throughout this section let $V'$ denote a vector space over
$\mathbb{K}$ such that $\dim V'=\dim V$.

\medskip

By a {\em $\mathbb{K}$-algebra isomorphism} from $\text{End}(V)$
to $\text{End}(V')$ we mean an isomorphism of $\mathbb{K}$-vector spaces
$\sigma : \text{End}(V) \to \text{End}(V')$
such that
$(XY)^\sigma=X^\sigma Y^\sigma$ for all $X,Y \in \text{End}(V)$.

\medskip

Let $\Phi=(A;\{E_i\}_{i=0}^d;A^*;\{E^*_i\}_{i=0}^d)$ denote
a tridiagonal system on $V$.
For a map $\sigma : \text{End}(V) \to \text{End}(V')$ we define
\[
 \Phi^\sigma := 
   (A^\sigma; \{ {E_i}^\sigma\}_{i=0}^d; 
    {A^*}^\sigma; \{{E^*_i}^\sigma\}_{i=0}^d).
\]

\medskip

\begin{definition}    \label{def:isom}   \samepage
Let $\Phi$ denote a tridiagonal system on $V$ and let
$\Phi'$ denote a tridiagonal system on $V'$.
By an {\em isomorphism of tridiagonal systems} from $\Phi$ to
$\Phi'$ we mean a $\mathbb{K}$-algebra isomorphism
$\sigma : \text{End}(V) \to \text{End}(V')$ such that
$\Phi^\sigma=\Phi'$.
We say $\Phi$ and $\Phi'$ are {\em isomorphic} whenever
there exists an isomorphism of tridiagonal systems from $\Phi$
to $\Phi'$.
\end{definition}

\medskip

It is useful to interpret the concept of isomorphism as follows.
Let $\gamma:V \to V'$ denote an isomorphism of $\mathbb{K}$-vector spaces.
Then there exists a $\mathbb{K}$-algebra isomorphism 
$\sigma :\text{End}(V) \to \text{End}(V')$ such that
$X^\sigma=\gamma X \gamma^{-1}$ for all $X \in \text{End}(V)$.
Conversely let $\sigma: \text{End}(V) \to \text{End}(V')$ denote a
$\mathbb{K}$-algebra isomorphism.
By the Skolem-Noether theorem \cite[Corollary 9.122]{Rot} there
exists an isomorphism of $\mathbb{K}$-vector spaces $\gamma: V \to V'$ such that
$X^\sigma = \gamma X \gamma^{-1}$ for all $X \in \text{End}(V)$.

\medskip

\begin{lemma}          \label{lem:isomalt}  \samepage
Let $\Phi=(A;\{E_i\}_{i=0}^d;A^*;\{E^*_i\}_{i=0}^d)$ denote a
tridiagonal system on $V$ and let
$\Phi'=(B;\{F_i\}_{i=0}^d; B^*;\{F^*_i\}_{i=0}^d)$ denote a
tridiagonal system on $V'$.
Then the following {\rm (i)}, {\rm (ii)} are equivalent.
\begin{itemize}
\item[\rm (i)]
$\Phi$ and $\Phi'$ are isomorphic.
\item[\rm (ii)]
There exists an isomorphism of $\mathbb{K}$-vector spaces
$\gamma : V \to V'$ such that
$\gamma A = B \gamma$, $\gamma A^* = B^* \gamma$,
and 
$\gamma E_i = F_i \gamma$, $\gamma E^*_i=F^*_i\gamma$ for
$0 \leq i \leq d$.
\end{itemize}
\end{lemma}

\begin{proof}
(i)$\Rightarrow$(ii):
Let $\sigma$ denote an isomorphism of tridiagonal systems
from $\Phi$ to $\Phi'$.
By our comments below Definition \ref{def:isom} 
there exists an isomorphism of $\mathbb{K}$-vector spaces 
$\gamma: V \to V'$ such that
$X^\sigma = \gamma X \gamma^{-1}$ for all $X \in \text{End}(V)$.
Observe that $\gamma$ satisfies the requirements of (ii).

(ii)$\Rightarrow$(i):
Define $\sigma : \text{End}(V) \to \text{End}(V')$ such that
$X^\sigma=\gamma X \gamma^{-1}$ for $X \in \text{End}(V)$.
Observe that $\sigma$ is an isormorphism from $\Phi$
to $\Phi'$.
\end{proof}

\section{The $D_4$ action}

\indent
Let $\Phi=(A; \{E_i\}_{i=0}^d; A^*; \{E^*_i\}_{i=0}^d)$
denote a tridigonal system on $V$.
Then each of the following is a tridiagonal system on $V$:
\begin{align*}
\Phi^{*}  &:= 
       (A^*; \{E^*_i\}_{i=0}^d; A; \{E_i\}_{i=0}^d), \\
\Phi^{\downarrow} &:=
       (A; \{E_i\}_{i=0}^d; A^*; \{E^*_{d-i}\}_{i=0}^d), \\
\Phi^{\Downarrow} &:=
       (A; \{E_{d-i}\}_{i=0}^d; A^*; \{E^*_{i}\}_{i=0}^d).
\end{align*}
Viewing $*$, $\downarrow$, $\Downarrow$ as permutations on the set of
all tridiagonal systems,
\begin{equation}    \label{eq:relation1}
\text{$*^2$$\;=\;$$\downarrow^2$$\;=\;$$\Downarrow^2$$\;=\;$$1$},
\end{equation}
\begin{equation}    \label{eq:relation2}
\text{$\Downarrow$$*$$\;=\;$$*$$\downarrow$}, \quad
\text{$\downarrow$$*$$\;=\;$$*$$\Downarrow$}, \quad
\text{$\downarrow$$\Downarrow$$\;=\;$$\Downarrow$$\downarrow$}.
\end{equation}
The group generated by symbols $*$, $\downarrow$, $\Downarrow$ subject
to the relations (\ref{eq:relation1}), (\ref{eq:relation2}) is the
dihedral group $D_4$. We recall that $D_4$ is the group of symmetries of a
square, and has $8$ elements.
Apparently $*$, $\downarrow$, $\Downarrow$ induce an action of $D_4$
on the set of all tridiagonal systems.
Two tridiagonal systems will be called {\em relatives} whenever they are
in the same orbit of this $D_4$ action. 
The relatives of $\Phi$ are as follows:

\medskip
\noindent
\begin{center}
\begin{tabular}{c|c}
name  &  relative \\
\hline
$\Phi$ & 
       $(A; \{E_i\}_{i=0}^d; A^*;  \{E^*_i\}_{i=0}^d)$ \\ 
$\Phi^{\downarrow}$ &
       $(A; \{E_i\}_{i=0}^d; A^*;  \{E^*_{d-i}\}_{i=0}^d)$ \\ 
$\Phi^{\Downarrow}$ &
       $(A; \{E_{d-i}\}_{i=0}^d; A^*;  \{E^*_i\}_{i=0}^d)$ \\ 
$\Phi^{\downarrow \Downarrow}$ &
       $(A; \{E_{d-i}\}_{i=0}^d; A^*;  \{E^*_{d-i}\}_{i=0}^d)$ \\ 
$\Phi^{*}$  & 
       $(A^*; \{E^*_i\}_{i=0}^d; A;  \{E_i\}_{i=0}^d)$ \\ 
$\Phi^{\downarrow *}$ &
       $(A^*; \{E^*_{d-i}\}_{i=0}^d; A;  \{E_i\}_{i=0}^d)$ \\ 
$\Phi^{\Downarrow *}$ &
       $(A^*; \{E^*_i\}_{i=0}^d; A;  \{E_{d-i}\}_{i=0}^d)$ \\ 
$\Phi^{\downarrow \Downarrow *}$ &
       $(A^*; \{E^*_{d-i}\}_{i=0}^d; A;  \{E_{d-i}\}_{i=0}^d)$
\end{tabular}
\end{center}

\medskip

From our comments below Note \ref{note:star} we obtain the following.

\medskip

\begin{lemma}              \label{lem:relshape}
Let $\Phi$ denote a tridiagonal system.
Then the following {\rm (i)}, {\rm (ii)} hold.
\begin{itemize}
\item[\rm (i)]
The relatives of $\Phi$ all have the same shape.
\item[\rm (ii)]
Suppose $\Phi$ is sharp.
Then each relative of $\Phi$ is sharp.
\end{itemize}
\end{lemma}

\medskip

We will use the following notational convention.

\medskip

\begin{definition}
Let $\Phi$ denote a tridiagonal system on $V$.
For $g \in D_4$ and for an object $f$ associated with $\Phi$ we let
$f^g$ denote the corresponding object associated with $\Phi^{g^{-1}}$.
\end{definition}

\section{The split decomposition}

\indent
In this section we recall the split decomposition associated with
a tridiagonal system \cite[Section 4]{ITT}.
With reference to Definition \ref{def}
for $0 \leq i \leq d$ we define
\begin{equation}          \label{eq:defUi}
   U_i = (E^*_0V+E^*_1V+\cdots+E^*_iV) \cap(E_iV+E_{i+1}V+\cdots+E_dV).
\end{equation}
By \cite[Theorem 4.6]{ITT}
\[
   V = U_0+U_1+\cdots+U_d  \qquad\qquad (\text{direct sum}),
\]
and for $0 \leq i \leq d$ both
\begin{align}
  U_0+U_1+\cdots+U_i &= E^*_0V+E^*_1V+\cdots+E^*_iV,    \label{eq:sumU0Ui} \\
  U_i+U_{i+1}+\cdots+U_d &= E_iV+E_{i+1}V+\cdots+E_dV.  \label{eq:sumUiUd}
\end{align}
By \cite[Corollary 5.7]{ITT}  $U_i$ has dimension $\rho_i$
for $0 \leq i \leq d$, where $\{\rho_i\}_{i=0}^d$ is the shape of $\Phi$.
By \cite[Theorem 4.6]{ITT} both
\begin{align}
  (A-\th_i I)U_i & \subseteq U_{i+1},    \label{eq:up}  \\
  (A^*-\th^*_i I)U_i & \subseteq U_{i-1}  \label{eq:down}
\end{align}
for $0 \leq i \leq d$, where $U_{-1}=0$ and $U_{d+1}=0$.
The sequence $\{U_i\}_{i=0}^d$ is called the {\em $\Phi$-split decomposition}
of $V$ \cite[Section 4]{ITT}.

\medskip

The following lemma will be useful.

\medskip

\begin{lemma}        \label{lem:bijection}       \samepage
With reference to Definition {\rm \ref{def}}
each of the following maps is bijective.
\begin{align*}
  E^*_0V &\to E_0V, \quad v \mapsto E_0v, &
  E^*_0V &\to E_dV, \quad v \mapsto E_dv,  \\
  E^*_dV &\to E_0V, \quad v \mapsto E_0v, &
  E^*_dV &\to E_dV, \quad v \mapsto E_dv,  \\
  E_0V &\to E^*_0V, \quad v \mapsto E^*_0v, &
  E_0V &\to E^*_dV, \quad v \mapsto E^*_dv, \\
  E_dV &\to E^*_0V, \quad v \mapsto E^*_0v, &
  E_dV &\to E^*_dV, \quad v \mapsto E^*_dv.
\end{align*}
\end{lemma}

\begin{proof}
Let $\xi$ denote the map on the right in the top line. 
We show that $\xi$ is bijective.
Let $\{U_i\}_{i=0}^d$ denote the $\Phi$-split decomposition of $V$.
By \cite[Lemmas 6.2, 6.5]{ITT}, the map $U_0 \to U_d$,
$u \mapsto \tau_d(A)u$ is a bijection. 
By \eqref{eq:defUi} we have $U_0=E^*_0V$ and $U_d=E_dV$. 
By \eqref{eq:E0} we have $\tau_d(A)=\tau_d(\th_d)E_d$.
By these comments $\xi$ is bijective. 
Applying $D_4$ we find each of the remaining maps is bijective.
\end{proof}

\section{Sharp tridiagonal systems and the parameter array}

\indent
Let $\Phi=(A;\{E_i\}_{i=0}^d;A^*;\{E^*_i\}_{i=0}^d)$ denote
a sharp tridiagonal system.
A bit later in the paper we will 
associate with $\Phi$ some formulae that have the scalars
\begin{equation}        \label{eq:traces}
 \tr(E_0E^*_0), \quad \tr(E_0E^*_d), \quad \tr(E_dE^*_0), \quad \tr(E_dE^*_d)
\end{equation}
in the denominator. So we take a moment to establish
that these scalars are nonzero.

\medskip

\begin{lemma}           \label{lem:trnonzero}      \samepage
With reference to Definition {\rm \ref{def}} assume $\Phi$ is sharp.
Then each of the traces \eqref{eq:traces} is nonzero.
\end{lemma}

\begin{proof}
We first show that $\tr(E_0E^*_0) \neq 0$.
Composing the first and third maps in the first column of
Lemma \ref{lem:bijection} we find that
the map $E_0V \to E_0V$, $v \mapsto E_0E^*_0v$ is bijective.
By this and since $\dim E_0V=1$, we find $E_0E^*_0E_0$ is a nonzero
scalar multiple of $E_0$.
Now take the trace and use $\tr(E_0)=1$
to find $\tr(E_0E^*_0) \neq 0$ as desired.
Applying this to the relatives of $\Phi$
we find each of the remaining traces is nonzero.
\end{proof}

\medskip

With reference to Definition \ref{def}
assume $\Phi$ is sharp, and let
$\{U_i\}_{i=0}^d$ denote the $\Phi$-split decomposition of $V$.
Note that $U_0$ has dimension $1$.
For $0 \leq i \leq d$ the space $U_0$ is invariant under
\begin{equation}        \label{eq:updown}
 (A^*-\th^*_1I)(A^*-\th^*_2I)\cdots(A^*-\th^*_iI)
 (A-\th_{i-1}I)\cdots(A-\th_1I)(A-\th_0I);
\end{equation}
let $\zeta_i$ denote the corresponding eigenvalue.
Note that $\zeta_0=1$.
We call the sequence $\{\zeta_i\}_{i=0}^d$ the {\em split sequence}
of $\Phi$.

\medskip

\begin{definition}  \samepage
With reference to Definition \ref{def}
assume $\Phi$ is sharp.
By the {\em parameter array} of $\Phi$ we mean the sequence
  $(\{\th_i\}_{i=0}^d; \{\th^*_i\}_{i=0}^d; \{\zeta_i\}_{i=0}^d)$
where $\{\zeta_i\}_{i=0}^d$ is the split sequence of $\Phi$.
We remark that the parameter array of $\Phi$ is defined only 
when $\Phi$ is sharp.
\end{definition}

\medskip

We now state a conjecture which indicates the
significance of the parameter array.

\medskip

\begin{conjecture}{\rm \cite{IT:Krawt}}   \label{conj:main}     \samepage
Let $d$ denote a nonnegative integer and let
\begin{equation}         \label{eq:parray}
 (\{\th_i\}_{i=0}^d; \{\th^*_i\}_{i=0}^d; \{\zeta_i\}_{i=0}^d)
\end{equation}
denote a sequence of scalars taken from $\mathbb{K}$.
Then there exists a sharp tridiagonal system $\Phi$ over $\mathbb{K}$
with parameter array \eqref{eq:parray} if and only if
{\rm (i)}--{\rm (iii)} hold below.
\begin{itemize}
\item[\rm (i)]
$\th_i \neq \th_j$, $\th^*_i \neq \th^*_j$ if $i \neq j$ $(0 \leq i,j \leq d)$.
\item[\rm (ii)]
$\zeta_0=1$, $\zeta_d \neq 0$, and
\begin{equation}           \label{eq:ineq}
 \sum_{i=0}^d \eta_{d-i}(\th_0)\eta^*_{d-i}(\th^*_0) \zeta_i \neq 0.
\end{equation}
\item[\rm (iii)]
The expressions
\[
  \frac{\th_{i-2}-\th_{i+1}}{\th_{i-1}-\th_i},  \qquad\qquad
  \frac{\th^*_{i-2}-\th^*_{i+1}}{\th^*_{i-1}-\th^*_i}
\]
are equal and independent of $i$ for $2 \leq i \leq d-1$.
\end{itemize}
Suppose {\rm (i)}--{\rm (iii)} hold. Then $\Phi$ is unique up to isomorphism of
tridiagonal systems.
\end{conjecture}

\medskip

Later in the paper we will prove one direction 
of Conjecture \ref{conj:main}.

\section{Some formulae for the split sequence}

\indent
In this section we obtain some formulae for the split sequence
of a sharp tridiagonal system.
To prepare for this we have a few lemmas which hold for a general
tridiagonal system.

\medskip

\begin{lemma}      \label{lem:E0tausiAstaujAEs0}            \samepage
With reference to Definition {\rm \ref{def}}, for $0 \leq i,j \leq d$
both
\begin{align}
 E_0 \tau^*_i(A^*)\tau_j(A)E^*_0 &= 0,      \label{eq:E0tausiAstaujAEs0} \\
 E^*_0 \tau_i(A)\tau^*_j(A^*)E_0 &= 0      \label{eq:Es0tauiAtausjAsE0}
\end{align}
provided $i \neq j$.
\end{lemma}

\begin{proof}
We first show \eqref{eq:E0tausiAstaujAEs0} for $i<j$.
In the left-hand side of \eqref{eq:E0tausiAstaujAEs0}
we insert a factor $I$ between $\tau^*_i(A^*)$ and $\tau_j(A)$.
We expand using $I=\sum_{r=0}^d E_r$ and simplify the result
using $E_r\tau_j(A)=\tau_j(\th_r)E_r$ for $0 \leq r \leq d$.
By these comments the left-hand side of \eqref{eq:E0tausiAstaujAEs0}
is equal to
\begin{equation}          \label{eq:sum1}
  \sum_{r=0}^d \tau_j(\th_r) E_0\tau^*_i(A^*)E_rE^*_0.
\end{equation}
For $0 \leq r \leq d$ we examine term $r$ in \eqref{eq:sum1}.
By the definition of $\tau_j$ we have $\tau_j(\th_r)=0$
for $0 \leq r \leq j-1$.
Using Lemma \ref{lem:trid}(i)
we find $E_0A^{*s}E_r=0$ for $0 \leq s \leq r-1$.
By this and since $\tau^*_i(A^*)$ has degree $i$ we find
$E_0\tau^*_i(A^*)E_r=0$ for $i+1 \leq r \leq d$.
By these comments term $r$ in \eqref{eq:sum1} vanishes for $0 \leq r \leq j-1$
and $i+1 \leq r \leq d$.
Recall $i<j$ so term $r$ in \eqref{eq:sum1} vanishes for $0 \leq r \leq d$.
In other words \eqref{eq:sum1} is $0$.
We have shown \eqref{eq:E0tausiAstaujAEs0} for $i<j$.
Next we show \eqref{eq:E0tausiAstaujAEs0} for $i>j$.
In the left-hand side of \eqref{eq:E0tausiAstaujAEs0}
we again insert a factor $I$ between $\tau^*_i(A^*)$ and $\tau_j(A)$.
This time we expand using $I=\sum_{r=0}^d E^*_r$ and simplify the result
using $\tau^*_i(A^*)E^*_r=\tau^*_i(\th^*_r)E^*_r$ for $0 \leq r \leq d$.
By these comments the left-hand side of \eqref{eq:E0tausiAstaujAEs0}
is equal to
\begin{equation}          \label{eq:sum2}
  \sum_{r=0}^d \tau^*_i(\th^*_r) E_0E^*_r\tau_j(A)E^*_0.
\end{equation}
For $0 \leq r \leq d$ we examine term $r$ in \eqref{eq:sum2}.
By the definition of $\tau^*_i$ we have $\tau^*_i(\th^*_r)=0$
for $0 \leq r \leq i-1$.
Using Lemma \ref{lem:trid}(ii)
we find $E^*_rA^s E^*_0=0$ for $0 \leq s \leq r-1$.
By this and since $\tau_j(A)$ has degree $j$ we find
$E^*_r \tau_j(A)E^*_0=0$ for $j+1 \leq r \leq d$.
By these comments term $r$ in \eqref{eq:sum2} vanishes for $0 \leq r \leq i-1$
and $j+1 \leq r \leq d$.
Recall $i>j$ so term $r$ in \eqref{eq:sum2} vanishes for $0 \leq r \leq d$.
In other words \eqref{eq:sum2} is $0$.
We have shown \eqref{eq:E0tausiAstaujAEs0} for $i>j$.
To get \eqref{eq:Es0tauiAtausjAsE0} apply \eqref{eq:E0tausiAstaujAEs0}
to $\Phi^*$.
\end{proof}

\begin{lemma}     \label{lem:formula1}   \samepage
With reference to Definition {\rm \ref{def}}
the following {\rm (i)}, {\rm (ii)} hold.
\begin{itemize}
\item[\rm (i)]
For $0 \leq i \leq d$ both
\begin{align}
E_0\tau^*_i(A^*)\tau_i(A)E^*_0 &=
 (\th_0-\th_1)(\th_0-\th_2)\cdots(\th_0-\th_i)
   E_0\tau^*_i(A^*)E_0E^*_0,            \label{eq:E0tausiAstauiAEs01}  \\
E_0\tau^*_i(A^*)\tau_i(A)E^*_0 &=
 (\th^*_0-\th^*_1)(\th^*_0-\th^*_2)\cdots(\th^*_0-\th^*_i)
   E_0E^*_0 \tau_i(A)E^*_0.             \label{eq:E0tausiAstauiAEs02}
\end{align}
\item[\rm (ii)]
For $0 \leq i \leq d$ both
\begin{align}
 E^*_0\tau_i(A)\tau^*_i(A^*)E_0 &=
   (\th^*_0-\th^*_1)(\th^*_0-\th^*_2)\cdots(\th^*_0-\th^*_i)
      E^*_0\tau_i(A)E^*_0E_0,           \label{eq:Es0tauiAtausiAsE01}  \\
 E^*_0\tau_i(A)\tau^*_i(A^*)E_0 &=
   (\th_0-\th_1)(\th_0-\th_2)\cdots(\th_0-\th_i)
     E^*_0E_0\tau^*_i(A^*)E_0.           \label{eq:Es0tauiAtausiAsE02} 
\end{align}
\end{itemize}
\end{lemma}

\begin{proof}
In the expression on the right in \eqref{eq:E0tausiAstauiAEs01},
eliminate the middle $E_0$ using the equation on the left in
\eqref{eq:E0};
evaluate the result using \eqref{eq:etad} and then
\eqref{eq:E0tausiAstaujAEs0} to get the expression on the left in
\eqref{eq:E0tausiAstauiAEs01}.
This gives \eqref{eq:E0tausiAstauiAEs01}.
Lines \eqref{eq:E0tausiAstauiAEs02}--\eqref{eq:Es0tauiAtausiAsE02}
are similarly obtained.
\end{proof}

\medskip

We now restrict our attention to sharp tridiagonal systems.

\medskip

\begin{theorem}           \label{thm:aux}   \samepage
With reference to Definition {\rm \ref{def}}, assume $\Phi$ is sharp and 
let $\{\zeta_i\}_{i=0}^d$ denote the split sequence of $\Phi$.
Then for $0 \leq i \leq d$ both
\begin{align} 
 (A^*-\th^*_1I)(A^*-\th^*_2I)\cdots(A^*-\th^*_iI)\tau_i(A)E^*_0
    &= \zeta_i E^*_0,                         \label{eq:auxAA}  \\
 (A-\th_1I)(A-\th_2I)\cdots(A-\th_iI)\tau^*_i(A^*)E_0
    &= \zeta_i E_0.                           \label{eq:auxAA2} 
\end{align}
Moreover $\zeta_i=\zeta^*_i$.
\end{theorem}

\begin{proof}
By the construction
\[
  (A^*-\th^*_1I)(A^*-\th^*_2I)\cdots(A^*-\th^*_iI)\tau_i(A)-\zeta_iI
\]
vanishes on $E^*_0V$, so \eqref{eq:auxAA} holds.
Next we show $\zeta_i=\zeta^*_i$.
In \eqref{eq:auxAA} take the trace of both sides and use 
$\tr(MN)=\tr(NM)$, $\tr(E^*_0)=1$, $E^*_0A^*=\th^*_0E^*_0$ to find
\begin{equation}                    \label{eq:auxzetai}
 (\th^*_0-\th^*_1)(\th^*_0-\th^*_2)\cdots(\th^*_0-\th^*_i) \tr(\tau_i(A)E^*_0)
  = \zeta_i.
\end{equation}
Applying this to $\Phi^*$,
\begin{equation}                   \label{eq:auxzetasi}
(\th_0-\th_1)(\th_0-\th_2)\cdots(\th_0-\th_i) \tr(\tau^*_i(A^*)E_0)
  = \zeta^*_i.
\end{equation}
To get $\zeta_i=\zeta^*_i$
we show the left-hand sides of \eqref{eq:auxzetai} and \eqref{eq:auxzetasi}
coincide.
Observe that the left-hand sides of \eqref{eq:E0tausiAstauiAEs01} and
\eqref{eq:E0tausiAstauiAEs02} coincide.
Since $E_0$ has rank $1$ we find $E_0E^*_0E_0$ is a scalar multiple of $E_0$.
Taking the trace we find $E_0E^*_0E_0=\tr(E_0E^*_0)E_0$. 
Using this and $\tr(MN)=\tr(NM)$ we find
that the trace of the right-hand side of \eqref{eq:E0tausiAstauiAEs01}
is equal to the left-hand side of \eqref{eq:auxzetasi} times $\tr(E_0E^*_0)$.
Similarly the trace of the right-hand side of \eqref{eq:E0tausiAstauiAEs02}
is equal to the left-hand side of \eqref{eq:auxzetai} times $\tr(E_0E^*_0)$.
By these comments the left-hand sides of
\eqref{eq:auxzetai}, \eqref{eq:auxzetasi} coincide so $\zeta_i=\zeta^*_i$.
To get \eqref{eq:auxAA2} apply \eqref{eq:auxAA} to $\Phi^*$ and use 
$\zeta_i=\zeta^*_i$.
\end{proof}

\begin{theorem}             \label{thm:E0tausiAstauiAEs0}
With reference to Definition {\rm \ref{def}}, assume $\Phi$ is sharp and 
let $\{\zeta_i\}_{i=0}^d$ denote the split sequence of $\Phi$.
Then for $0 \leq i \leq d$ both
\begin{align} 
 E_0\tau^*_i(A^*)\tau_i(A)E^*_0 &= \zeta_i E_0E^*_0, 
                                      \label{eq:E0tausiAstauiAEs0} \\
 E^*_0\tau_i(A)\tau^*_i(A^*)E_0 &= \zeta_i E^*_0E_0. 
                                      \label{eq:Es0tauiAtausiAsE0}
\end{align}
\end{theorem}

\begin{proof}
We first show \eqref{eq:E0tausiAstauiAEs0}.
Multiplying both sides of \eqref{eq:auxAA} on the left by $E_0$,
\begin{equation}             \label{eq:auxA}
 E_0(A^*-\th^*_1I)(A^*-\th^*_2I)\cdots(A^*-\th^*_iI)\tau_i(A)E^*_0
   =\zeta_iE_0E^*_0.
\end{equation}
Observe that the expression
\[
 (A^*-\th^*_1I)(A^*-\th^*_2I)\cdots(A^*-\th^*_iI) - \tau^*_i(A^*)
\]
is a polynomial in $A^*$ with degree less than $i$,
so it is a linear combination of $\{\tau^*_r(A^*)\}_{r=0}^{i-1}$.
By this and Lemma \ref{lem:E0tausiAstaujAEs0}
we find that the left-hand side of \eqref{eq:auxA} is equal
to the left-hand side of \eqref{eq:E0tausiAstauiAEs0}.
This gives \eqref{eq:E0tausiAstauiAEs0}.
To get \eqref{eq:Es0tauiAtausiAsE0} apply \eqref{eq:E0tausiAstauiAEs0}
to $\Phi^*$ and use $\zeta_i=\zeta^*_i$.
\end{proof}

\medskip

In the proof of Theorem \ref{thm:aux} we used some trace formula
for $\zeta_i$.
Taking the trace in \eqref{eq:E0tausiAstauiAEs0}, 
\eqref{eq:Es0tauiAtausiAsE0} we obtain some more trace formulae for $\zeta_i$.
These formulae are summarized below.

\medskip

\begin{theorem}              \label{thm:zetai}  \samepage
With reference to Definition {\rm \ref{def}}, assume $\Phi$ is sharp and 
let $\{\zeta_i\}_{i=0}^d$ denote the split sequence of $\Phi$.
Then {\rm (i)}, {\rm (ii)} hold below.
\begin{itemize}
\item[\rm (i)]
For $0 \leq i \leq d$,
\begin{align}
 \zeta_i &= (\th^*_0-\th^*_1)(\th^*_0-\th^*_2)\cdots(\th^*_0-\th^*_i)
             \tr(\tau_i(A)E^*_0),                   \label{eq:zetai1}  \\
 \zeta_i &= (\th_0-\th_1)(\th_0-\th_2)\cdots(\th_0-\th_i)
            \tr(\tau^*_i(A^*)E_0).                  \label{eq:zetai4}
\end{align}
\item[\rm (ii)]
For $0 \leq i \leq d$,
\begin{align}
 \zeta_i &= \frac{\tr(E_0\tau^*_i(A^*)\tau_i(A)E^*_0)}
                  {\tr(E_0E^*_0)},            \label{eq:zetai2}  \\
 \zeta_i &= \frac{\tr(E^*_0\tau_i(A)\tau^*_i(A)E_0)}
                 {\tr(E^*_0E_0)}.                 \label{eq:zetai3}
\end{align}
\end{itemize}
\end{theorem}

\section{The split sequence for the relatives of $\Phi$, part I}

\indent
We now discuss the relationship between the split sequences
for the relatives of $\Phi$. 
In this discussion we treat separately the last term
in the split sequence, since its role is somewhat special
as we shall see.
In this section we treat the last term.
In Section \ref{sec:zetai} we
treat the remaining terms.
In this section we also give a proof of Conjecture \ref{conj:main} in
one direction.

\medskip

\begin{theorem}       \label{thm:zetad}     \samepage
With reference to Definition {\rm \ref{def}}, assume $\Phi$ is sharp and let
$\{\zeta_i\}_{i=0}^d$ denote the split sequence of $\Phi$.
Then both
\begin{equation}         \label{eq:zetad}
\zeta_d= \eta^*_d(\th^*_0)\tau_d(\th_d)\tr(E_dE^*_0),
\qquad\qquad
\zeta_d= \eta_d(\th_0)\tau^*_d(\th^*_d)\tr(E^*_dE_0).
\end{equation}
\end{theorem}

\begin{proof}
To get the equation on the left in \eqref{eq:zetad}, set $i=d$ in
\eqref{eq:zetai1} and evaluate the result using the equation on the right in
\eqref{eq:E0}. 
The equation on the right in \eqref{eq:zetad} is similarly obtained
using \eqref{eq:zetai4}.
\end{proof}

\begin{theorem}     \label{thm:zetadrel}    \samepage
With reference to Definition {\rm \ref{def}}, assume $\Phi$ is sharp and let
$\{\zeta_i\}_{i=0}^d$ denote the split sequence of $\Phi$.
Then {\rm (i)}, {\rm (ii)} hold below.
\begin{itemize}
\item[\rm (i)]
For the tridiagonal systems 
$\Phi$, $\Phi^*$, $\Phi^{\d\D}$, $\Phi^{\d\D*}$
the last term in the split sequence is equal to $\zeta_d$.
\item[\rm (ii)]
For the tridiagonal systems
$\Phi^\d$, $\Phi^\D$, $\Phi^{\d*}$, $\Phi^{\D*}$
the last term in the split sequence is equal to
\begin{equation}         \label{eq:zetadD}
\sum_{i=0}^d \eta_{d-i}(\th_0)\eta^*_{d-i}(\th^*_0)\zeta_i.
\end{equation}
\end{itemize}
\end{theorem}

\begin{proof}
(i):
We have $\zeta^*_d=\zeta_d$ by Theorem \ref{thm:aux}.
Applying the equation on the left in \eqref{eq:zetad} to $\Phi^{\d\D}$
we find
\[
   {\zeta_d}^{\d\D}= \tau^*_d(\th^*_d)\eta_d(\th_0)\tr(E_0E^*_d).
\]
Comparing this with the equation on the right in \eqref{eq:zetad} we get 
${\zeta_d}^{\d\D}=\zeta_d$.
By these comments we get
$\zeta_d={\zeta_d}^*={\zeta_d}^{\d\D}={\zeta_d}^{\d\D*}$.

(ii):
Applying (i) to $\Phi^{\D}$ we get
${\zeta_d}^{\D}={\zeta_d}^{\D*}={\zeta_d}^{\d}={\zeta_d}^{\d*}$,
so it suffices to show that the expression \eqref{eq:zetadD} is equal to
${\zeta_d}^\D$.
By \eqref{eq:zetai1} at $i=d$,
\[
  \zeta_d=\eta^*_d(\th^*_0)\tr(\tau_d(A)E^*_0).
\]
Applying this to $\Phi^\D$,
\[
 {\zeta_d}^{\D} = \eta^*_d(\th^*_0)\tr(\eta_d(A)E^*_0).
\]
Evaluating \eqref{eq:etad} at $\lambda=A$,
\[
  \eta_d(A)=\sum_{i=0}^d \eta_{d-i}(\th_0)\tau_i(A).
\]
By these comments and \eqref{eq:zetai1},
\begin{align*}
 {\zeta_d}^{\D}
 &= \eta^*_d(\th^*_0)\sum_{i=0}^d \eta_{d-i}(\th_0)\tr(\tau_i(A)E^*_0) \\
 &= \sum_{i=0}^d \eta_{d-i}(\th_0)\eta^*_{d-i}(\th^*_0)
          (\th^*_0-\th^*_1)\cdots(\th^*_0-\th^*_i)\tr(\tau_i(A)E^*_0) \\
 &= \sum_{i=0}^d \eta_{d-i}(\th_0)\eta^*_{d-i}(\th^*_0)\zeta_i.
\end{align*}
So (ii) holds.
\end{proof}

\begin{corollary}        \label{cor:zetad}      \samepage
With reference to Definition {\rm \ref{def}}, assume $\Phi$ is sharp and let
$\{\zeta_i\}_{i=0}^d$ denote the split sequence of $\Phi$.
Then $\zeta_0=1$, $\zeta_d \neq 0$, and \eqref{eq:ineq} holds.
\end{corollary}

\begin{proof}
Obviously $\zeta_0=1$.
By  Lemma \ref{lem:trnonzero} and Theorem \ref{thm:zetad} we get
$\zeta_d \neq 0$.
Applying this to $\Phi^\D$ we find ${\zeta_d}^\D \neq 0$.
Now \eqref{eq:ineq} follows from this and Theorem \ref{thm:zetadrel}(ii).
\end{proof}

\medskip

We can now easily prove Conjecture \ref{conj:main} in one direction.

\medskip

\begin{proofof}{Conjecture \ref{conj:main} (direction ``only if'')}
Assume that there exists a sharp tridiagonal system $\Phi$ over $\mathbb{K}$
that has parameter array \eqref{eq:parray}. 
We show that this parameter array satisfies (i)--(iii).
Assertion (i) follows from Definition \ref{def}.
Assertion (ii) is just Corollary \ref{cor:zetad}.
Assertion (iii) is just Lemma \ref{lem:indep}.
\end{proofof}

\section{The split sequence for the relatives of $\Phi$, part II}
\label{sec:zetai}

\indent
In this section we continue to discuss the relationship
between the split sequences for the relatives of $\Phi$.
We will need the following scalars.
Given a tridiagonal system
$(A;\{E_i\}_{i=0}^d;A^*;\{E^*_i\}_{i=0}^d)$ over $\mathbb{K}$
and given nonnegative integers $r,s,t$
such that $r+s+t \leq d$,
in \cite[Definition 13.1]{T:24points} we defined a scalar
$[r,s,t]_q \in \mathbb{K}$, 
where $q+q^{-1}+1$ is the common value of \eqref{eq:beta}.
For example, if $q \neq 1$ and $q \neq -1$ then
\[
  [r,s,t]_q =
  \frac{(q;q)_{r+s}(q;q)_{r+t}(q;q)_{s+t}}
       {(q;q)_r (q;q)_s (q;q)_t(q;q)_{r+s+t}},
\]
where
\[
  (a;q)_n=(1-a)(1-aq)\cdots(1-aq^{n-1}).
\]

\medskip

\begin{lemma} {\rm \cite[Theorem 5.5]{NT:formula}} \label{lem:formula}  \samepage
With reference to Definition {\rm \ref{def}},
for $0 \leq i \leq d$ we have
\begin{equation}         \label{eq:etaitauh}
 \eta_i = \sum_{h=0}^i [h,i-h,d-i]_q \eta_{i-h}(\th_0)\tau_h.
\end{equation}
\end{lemma}

\begin{note}
In \cite{NT:formula} we gave a proof of \eqref{eq:etaitauh}  for the
case of $q \neq 1$, $q \neq -1$. 
A similar proof establishes \eqref{eq:etaitauh} for each of
the following cases:
(i) $q \neq 1$, $q \neq -1$;
(ii) $q =1$, $\text{Char}(\mathbb{K}) \neq 2$;
(iii) $q=-1$, $\text{Char}(\mathbb{K}) \neq 2$;
(iv) $q=1$, $\text{Char}(\mathbb{K})=2$.
\end{note}

\medskip

We are now ready to give the relationship between the
split sequences for the relatives of $\Phi$.
These relationships will show that the split sequence
for each relative of $\Phi$ is determined by the parameter
array of $\Phi$.
We start with a comment. By the last line of Theorem \ref{thm:aux},
in each column of the following array
the relatives of $\Phi$ have the same split sequence:
\[
\begin{array}{llll}
 \Phi & \Phi^\d & \Phi^\D & \Phi^{\d\D}  \\ 
 \Phi^* & \Phi^{\d *} & \Phi^{\D *} & \Phi^{\d\D *}
\end{array}
\]
Therefore we limit our attention to the split sequences for the
relatives of $\Phi$ in the first row.

\medskip

\begin{theorem}           \samepage
With reference to Definition {\rm \ref{def}},
assume $\Phi$ is sharp and let $\{\zeta_i\}_{i=0}^d$ denote the
split sequence of $\Phi$.
Then {\rm (i)}--{\rm (iv)} hold below.
\begin{itemize}
\item[\rm (i)]
For $0 \leq i \leq d$ both
\begin{align*}
 \frac{\zeta_i^{\d}}
      {(\th_0-\th_{1})(\th_0-\th_{2})\cdots(\th_0-\th_{i})}
 &= \sum_{h=0}^i
  \frac{[h,i-h,d-i]_q\eta^*_{i-h}(\th^*_0)\zeta_h}
       {(\th_0-\th_{1})(\th_0-\th_{2})\cdots(\th_0-\th_{h})},  \\
 \frac{\zeta_i}
      {(\th_0-\th_{1})(\th_0-\th_{2})\cdots(\th_0-\th_{i})}
 &= \sum_{h=0}^i
  \frac{[h,i-h,d-i]_q\tau^*_{i-h}(\th^*_d)\zeta_h^{\d}}
       {(\th_0-\th_{1})(\th_0-\th_{2})\cdots(\th_0-\th_{h})}.
\end{align*}
\item[\rm (ii)]
For $0 \leq i \leq d$ both
\begin{align*}
 \frac{\zeta_i^{\D}}
      {(\th^*_0-\th^*_{1})(\th^*_0-\th^*_{2})\cdots(\th^*_0-\th^*_{i})}
 &= \sum_{h=0}^i
  \frac{[h,i-h,d-i]_q\eta_{i-h}(\th_0)\zeta_h}
       {(\th^*_0-\th^*_{1})(\th^*_0-\th^*_{2})\cdots(\th^*_0-\th^*_{h})},  \\
 \frac{\zeta_i}
      {(\th^*_0-\th^*_{1})(\th^*_0-\th^*_{2})\cdots(\th^*_0-\th^*_{i})}
 &= \sum_{h=0}^i
  \frac{[h,i-h,d-i]_q\tau_{i-h}(\th_d)\zeta_h^{\D}}
       {(\th^*_0-\th^*_{1})(\th^*_0-\th^*_{2})\cdots(\th^*_0-\th^*_{h})}.
\end{align*}
\item[\rm (iii)]
For $0 \leq i \leq d$ both
\begin{align*}
 \frac{\zeta_i^{\d\D}}
      {(\th^*_d-\th^*_{d-1})(\th^*_d-\th^*_{d-2})\cdots(\th^*_d-\th^*_{d-i})}
 &= \sum_{h=0}^i
 \frac{[h,i-h,d-i]_q\eta_{i-h}(\th_0)\zeta_h^{\d}}
      {(\th^*_d-\th^*_{d-1})(\th^*_d-\th^*_{d-2})\cdots(\th^*_d-\th^*_{d-h})},\\
 \frac{\zeta_i^{\d}}
      {(\th^*_d-\th^*_{d-1})(\th^*_d-\th^*_{d-2})\cdots(\th^*_d-\th^*_{d-i})}
 &= \sum_{h=0}^i
  \frac{[h,i-h,d-i]_q\tau_{i-h}(\th_d)\zeta_h^{\d\D}}
       {(\th^*_d-\th^*_{d-1})(\th^*_d-\th^*_{d-2})\cdots(\th^*_d-\th^*_{d-h})}.
\end{align*}
\item[\rm (iv)]
For $0 \leq i \leq d$ both
\begin{align*}
 \frac{\zeta_i^{\d\D}}
      {(\th_d-\th_{d-1})(\th_d-\th_{d-2})\cdots(\th_d-\th_{d-i})}
 &= \sum_{h=0}^i
  \frac{[h,i-h,d-i]_q\eta^*_{i-h}(\th^*_0)\zeta_h^{\D}}
       {(\th_d-\th_{d-1})(\th_d-\th_{d-2})\cdots(\th_d-\th_{d-h})}, \\
 \frac{\zeta_i^{\D}}
      {(\th_d-\th_{d-1})(\th_d-\th_{d-2})\cdots(\th_d-\th_{d-i})}
 &= \sum_{h=0}^i
  \frac{[h,i-h,d-i]_q\tau^*_{i-h}(\th^*_d)\zeta_h^{\d\D}}
       {(\th_d-\th_{d-1})(\th_d-\th_{d-2})\cdots(\th_d-\th_{d-h})}.
\end{align*}
\end{itemize}
\end{theorem}

\begin{proof}
We start by obtaining the first equation in part (ii).
Applying \eqref{eq:zetai1} to $\Phi^\D$,
\begin{equation}         \label{eq:zetaaux1}
 \zeta_i^{\D} = 
  (\th^*_0-\th^*_1)(\th^*_0-\th^*_2)\cdots(\th^*_0-\th^*_i)
   \tr(\eta_i(A)E^*_0).
\end{equation}
Evaluate \eqref{eq:etaitauh} at $\lambda=A$ and in the result multiply 
both sides on the right by $E^*_0$ to get
\[
 \eta_i(A)E^*_0 = \sum_{h=0}^i [h,i-h,d-i]_q\eta_{i-h}(\th_0)\tau_h(A)E^*_0.
\]
In this equation we take the trace of both sides and use
\eqref{eq:zetai1}, \eqref{eq:zetaaux1} to obtain the first equation
in part (ii).
Apply $D_4$ to this and use $\zeta_i=\zeta^*_i$
to obtain the remaining formulae.
\end{proof}

\section{Bilinear forms, anti-automorphisms, and tridiagonal systems}

\indent
With reference to Definition \ref{def}, assume for the moment that
$\Phi$ is sharp.
Our next goal is to show that if Conjecture \ref{conj:main} is true
then there exists a nondegenerate symmetric bilinear form
$\b{\;,\,}$ on $V$ that satisfies
\begin{equation}            \label{eq:compatible}
  \b{Au,v}=\b{u,Av},  \qquad
  \b{A^*u,v}=\b{u,A^*v}  \qquad\qquad  \text{for all $u,v \in V$}.
\end{equation}
We will also obtain some related results involving anti-automorphisms. 
We start with some definitions.
Throughout this section let $V'$ denote a vector space 
over $\mathbb{K}$ such that $\dim V'=\dim V$.

A map $\b{\;,\,} : V \times V' \to \mathbb{K}$ is called a
{\em bilinear form} whenever the following conditions hold
for $u,v \in V$, for $u',v' \in V'$, and for $\alpha \in \mathbb{K}$:
(i) $\b{u+v,u'}=\b{u,u'}+\b{v,u'}$;
(ii) $\b{\alpha u,u'}=\alpha \b{u,u'}$;
(iii) $\b{u,u'+v'}=\b{u,u'}+\b{u,v'}$;
(iv) $\b{u,\alpha u'}=\alpha \b{u,u'}$.
We observe that a scalar multiple of a bilinear form is a bilinear form.
Let $\b{\;,\,} : V \times V' \to \mathbb{K}$ denote a bilinear form.
Then the following are equivalent:
(i) there exists a nonzero $v \in V$ such that $\b{v,v'}=0$ for all $v'\in V'$;
(ii) there exists a nonzero $v' \in V'$ such that $\b{v,v'}=0$ for all $v\in V$.
The form $\b{\;,\,}$ is said to be {\em degenerate} whenever (i), (ii) hold
and {\em nondegenerate} otherwise.
By a {\em bilinear form on $V$} we mean a bilinear form
$\b{\;,\,}:V \times V \to \mathbb{K}$.
This form is said to be {\em symmetric} whenever $\b{u,v}=\b{v,u}$
for all $u,v \in V$.

\medskip

\begin{lemma}            \label{lem:nondeg}  \samepage
With reference to Definition {\rm \ref{def}},
let $\b{\;,\,}$ denote a nonzero bilinear form on $V$ that satisfies
\eqref{eq:compatible}.
Then $\b{\;,\,}$ is nondegenerate.
\end{lemma}

\begin{proof}
It suffices to show that the space
$W=\{w \in V \,|\, \b{w,V}=0\}$ is zero.
Using \eqref{eq:compatible} we routinely find
$AW \subseteq W$ and $A^*W \subseteq W$, so either
$W=0$ or $W=V$ by Definition \ref{def:TDpair}(iv).
But $W \neq V$ since  $\b{\;,\,}$ is nonzero
so $W=0$ as desired.
\end{proof}

\begin{lemma}          \label{lem:orth}  \samepage
With reference to Definition {\rm \ref{def}},
let $\b{\;,\,}$ denote a nonzero bilinear form on $V$ that satisfies 
\eqref{eq:compatible}.
Then {\rm (i)}, {\rm (ii)} hold below.
\begin{itemize}
\item[\rm (i)]
$\b{E_iV,E_jV}=0$ and $\b{E^*_iV,E^*_jV}=0$
if $i \neq j$ $(0 \leq i,j \leq d)$.
\item[\rm (ii)]
For $0 \leq i \leq d$
the restriction of $\b{\;,\,}$ to each of $E_iV$, $E^*_iV$ is nondegenerate.
\end{itemize}
\end{lemma}

\begin{proof}
(i):
Let $i,j$ be given with $i \neq j$. 
We show $\b{E_iu, E_jv}=0$ for $u,v\in V$. 
Recall that $E_i$ is contained in the subalgebra
of $\text{End}(V)$ generated by $A$, so using \eqref{eq:compatible} we have
$\b{E_iu,E_jv}=\b{u,E_iE_jv}=\b{u,0}=0$.
The proof of $\b{E^*_iV,E^*_jV}=0$ is similar.

(ii):
 Combine (i) above with Lemma \ref{lem:nondeg}.
\end{proof}

\begin{lemma}       \label{lem:unique}  \samepage
With reference to Definition {\rm \ref{def}} assume $\Phi$ is sharp.
Then up to scalar multiple there exists at most one bilinear form $\b{\;,\,}$
on $V$ that satisfies \eqref{eq:compatible}.
\end{lemma}

\begin{proof}
The dimension of $E_0V$ is $1$ since $\Phi$ is sharp; 
pick a nonzero $\eta \in E_0V$.
Let $\b{\;,\,}$ denote a nonzero bilinear form on $V$ that
satisfies \eqref{eq:compatible}, and note that 
$\b{\eta,\eta} \neq 0$ by Lemma \ref{lem:orth}(ii). 
Suppose another bilinear from $\b{\;,\,}'$ on $V$ 
satisfies \eqref{eq:compatible}. 
We show that $\b{\;,\,}'$ is a scalar multiple of $\b{\;,\,}$.
Since  $\b{\eta,\eta} \neq 0$ there exists $\alpha \in \mathbb{K}$
such that $\b{\eta,\eta}'=\alpha \b{\eta,\eta}$.
Define a map $(\;,\,) : V \times V \to \mathbb{K}$ by 
\[
 (u,v)= \b{u,v}'-\alpha \b{u,v}    \qquad \qquad  (u,v \in V).
\]
The map $(\;,\,)$ is a bilinear form on $V$ that satisfies 
\eqref{eq:compatible}
and $(\eta, \eta)=0$ so $(\;,\,)$ is zero by our preliminary comment.
Now $\b{\;,\,}'=\alpha \b{\;,\,}$ and the result follows.
\end{proof}

\begin{lemma}               \label{lem:symmetric}  \samepage
With reference to Definition {\rm \ref{def}} assume $\Phi$ is sharp.
Let $\b{\;,\,}$ denote a bilinear form on $V$ that satisfies
\eqref{eq:compatible}.
Then $\b{\;,\,}$ is symmetric.
\end{lemma}

\begin{proof}
We assume $\b{\;,\,}$ is nonzero; otherwise we are done.
Pick a nonzero $\eta \in E_0V$ and note that
$\b{\eta,\eta} \neq 0$ by Lemma \ref{lem:orth}(ii).
Define a map $(\;,\,): V \times V \to \mathbb{K}$ by
$(u,v)=\b{v,u}$ for $u,v \in V$. 
Then $(\;,\,)$ is a bilinear form on $V$ that satisfies \eqref{eq:compatible}, 
so by Lemma \ref{lem:unique} there exists $\alpha \in \mathbb{K}$ 
such that $(\;,\,)= \alpha  \b{\;,\,}$.
In other words $\b{u,v}=\alpha \b{v,u}$ for all $u,v \in V$.
Now $\b{\eta,\eta}=\alpha \b{\eta,\eta}$ and $\b{\eta,\eta} \neq 0$ so
$\alpha=1$.
Therefore $\b{u,v}=\b{v,u}$ for all $u,v \in V$,
so $\b{\;,\,}$ is symmetric.
\end{proof}

\medskip

By a {\em $\mathbb{K}$-algebra anti-isomorphism} from $\text{End}(V)$
to $\text{End}(V')$ we mean an isomorphism of $\mathbb{K}$-vector spaces
$\sigma : \text{End}(V) \to \text{End}(V')$ such that
$(XY)^\sigma=Y^\sigma X^\sigma$ for all $X,Y \in \text{End}(V)$.
By an {\em anti-automorphism} of $\text{End}(V)$ we mean
a $\mathbb{K}$-algebra anti-isomorphism from $\text{End}(V)$ to $\text{End}(V)$.
Bilinear forms and anti-isomorphisms are related as follows.
Let $\b{\;,\,} : V\times V' \to \mathbb{K}$ denote a nondegenerate
bilinear form. 
Then there exists a unique anti-isomorphism 
$\sigma :\text{End}(V) \to \text{End}(V')$ such that
$\b{Xv,v'}=\b{v,X^\sigma v'}$ for all $v \in V$, $v' \in V'$, 
$X \in \text{End}(V)$. 
Conversely, given an anti-isomorphism 
$\sigma :\text{End}(V) \to \text{End}(V')$ there exists a bilinear
form $\b{\;,\,} : V \times V' \to \mathbb{K}$ 
such that $\b{Xv,v'}=\b{v,X^\sigma v'}$ for all
$v \in V$, $v' \in V'$, $X \in \text{End}(V)$. This bilinear form is
nondegenerate, and 
uniquely determined by $\sigma$  up to multiplication by a nonzero
scalar in $\mathbb{K}$.
We say the form $\b{\;,\,}$ is {\em associated} with $\sigma$.

\medskip

\begin{lemma}              \label{lem:involutive}  \samepage
Let $\b{\;,\,}$ denote a nondegenerate bilinear form on $V$
and let $\sigma$ denote the anti-automorphism of $\text{\rm End}(V)$
associated with $\b{\;,\,}$.
Assume $\b{\;,\,}$ is symmetric.
Then $X^{\sigma\sigma} = X$ for all $X \in \text{\rm End}(V)$.
\end{lemma}

\begin{proof}
We fix $u \in V$ and show $(X^{\sigma \sigma}-X)u=0$.
For all $v \in V$ we have
$\b{Xu,v} 
 = \b{u,X^\sigma v} 
 = \b{X^\sigma v,u} 
 = \b{v, (X^\sigma)^\sigma u}
 = \b{X^{\sigma\sigma} u,v}$,
so $\b{(X-X^{\sigma\sigma})u,v}=0$.
Therefore $(X^{\sigma \sigma}-X)u=0$ 
since $\b{\;,\,}$ is nondegenerate.
\end{proof}

\begin{lemma}           \label{lem:antitr}  \samepage
Let $\sigma$ denote an anti-isomorphism from $\text{\rm End}(V)$
to $\text{\rm End}(V')$.
Then $\tr(X)=\tr(X^\sigma)$ for all $X \in \text{\rm End}(V)$.
\end{lemma}

\begin{proof}
Fix a basis $\{v_i\}_{i=1}^n$ for $V$ and a basis $\{v'_i\}_{i=1}^n$
for $V'$. 
For the purpose of this proof we identify
each element of $\text{End}(V)$ (resp. $\text{End}(V')$) with
the matrix in $\Mat{n}$ that represents it with
respect to $\{v_i\}_{i=1}^n$ (resp. $\{v'_i\}_{i=1}^n$).
By the Skolem-Noether theorem \cite[Corollary 9.122]{Rot} there exists
an invertible $R \in \Mat{n}$ such that
$X^{\sigma} = R^{-1}X^t R$ for all $X \in \text{End}(V)$.
Now $\tr(X^\sigma)=\tr(R^{-1}X^tR)=\tr(X^t)=\tr(X)$
as desired.
\end{proof}

\medskip

The following is a mild generalization of
\cite[Theorem 1.2]{AC3}.

\medskip

\begin{proposition} {\rm \cite[Theorem 1.2]{AC3}} \label{prop:antiTD} \samepage
With reference to Definition {\rm \ref{def}}, let $V'$ denote a vector space
such that $\dim V'=\dim V$, and
let $\sigma$ denote an anti-isomorphism from
$\text{\rm End}(V)$ to $\text{\rm End}(V')$.
Then {\rm (i)}--{\rm (iii)} hold below.
\begin{itemize}
\item[\rm (i)]
$\Phi^\sigma$ is a tridiagonal system on $V'$.
\item[\rm (ii)]
$\Phi$ and $\Phi^\sigma$ have the same
eigenvalue sequence and dual eigenvalue sequence.
\item[\rm (iii)]
Assume $\Phi$ is sharp.
Then $\Phi^\sigma$ is sharp.
Moreover $\Phi$ and $\Phi^\sigma$ have the same split sequence.
\end{itemize}
\end{proposition}

\begin{proof}
(i):
We first show that the pair $A^\sigma,A^{*\sigma}$ is a tridiagonal pair
on $V'$.
To do this we verify that $A^\sigma, A^{*\sigma}$ satisfy conditions 
(i)--(iv) of Definition \ref{def:TDpair}.
Concerning Definition \ref{def:TDpair}(i), 
we claim that $A^\sigma$ is diagonalizable with eigenvalues 
$\{\th_i\}_{i=0}^d$.
For $f \in \mathbb{K}[\lambda]$ we have $f(A)=0$ if and only if
$f(A^\sigma)=0$. 
Therefore $A$ and $A^\sigma$ have the same minimal polynomial.
The minimal polynomial of $A$ is $\prod_{i=0}^d (\lambda - \th_i)$
so the minimal polynomial of $A^\sigma$ is
$\prod_{i=0}^d (\lambda - \th_i)$. 
By this and since $\{\th_i\}_{i=0}^d$ are mutually distinct,
$A^\sigma$ is diagonalizable with eigenvalues $\{\th_i\}_{i=0}^d$.
Similarly $A^{*\sigma}$ is diagonalizable with eigenvalues $\{\th^*_i\}_{i=0}^d$.
Concerning Definition \ref{def:TDpair}(ii), 
for $0 \leq i \leq d$ we apply $\sigma$ to \eqref{eq:defEi}
and find $E_i^\sigma$ is the primitive idempotent of $A^\sigma$ 
associated with $\th_i$.
Define $V_i = E_i^\sigma V'$ and note that $\{V_i\}_{i=0}^d$ is an ordering
of the eigenspaces of $A^\sigma$. 
Applying $\sigma$ to the equation in Lemma \ref{lem:trid}(i) we find
$E_j^\sigma A^{*\sigma} E_i^\sigma = 0$ if $|i-j|>1$ $(0 \leq i,j\leq d)$.
Therefore $A^{*\sigma}V_i \subseteq V_{i-1} + V_i + V_{i+1}$ for
$0 \leq i \leq d$, where $V_{-1}=0$ and $V_{d+1}=0$.
The verification of Definition \ref{def:TDpair}(iii) is similar,
using $V^*_i = E^{*\sigma}_i V'$ for $0 \leq i \leq d$.
Concerning Definition \ref{def:TDpair}(iv), let $W$ denote a subspace
of $V'$ such that $A^\sigma W \subseteq W$ and $A^{*\sigma} W \subseteq W$.
We show $W=0$ or $W=V'$. 
Let $\b{\;,\,} : V \times V' \to \mathbb{K}$
denote the bilinear form associated with $\sigma$.
Define 
$W^\perp= \{v \in V \,|\, \b{v,w}=0 \text{ for all } w \in W \}$.
Since $\b{\;,\,}$ is nondegenerate 
$\dim W + \dim W^\perp$ is equal to the common dimension of $V, V'$. 
Observe $AW^\perp \subseteq W^\perp$; indeed
for all $w \in W$ and all $v \in W^\perp$, 
$\b{Av,w} = \b{v,A^\sigma w} = 0$ since $A^\sigma w \in W$.
We similarly obtain $A^*W^\perp \subseteq W^\perp$. 
Now $W^\perp =0$ or $W^\perp = V$ since $A,A^*$ is a tridiagonal pair on $V$, 
and therefore $W=V'$ or $W=0$.
We have shown the pair $A^\sigma, A^{*\sigma}$ satisfies conditions
(i)--(iv) of Definition \ref{def:TDpair}, so $A^\sigma, A^{*\sigma}$ 
is a tridiagonal pair on $V'$. 
By the construction $\{E^\sigma_i\}_{i=0}^d$
(resp. $\{E^{*\sigma}_i\}_{i=0}^d$) is a standard ordering of the 
primitive idempotents of $A^\sigma$ (resp. $A^{*\sigma}$).
Now $(A^\sigma;\{E^\sigma_i\}_{i=0}^d;A^*;\{E^{*\sigma}_i\})$
is a tridiagonal system on $V'$ by Definition \ref{def:TDsystem}.

(ii):
We mentioned in the proof of (i) above that for $0 \leq i \leq d$
the scalar $\th_i$ is the eigenvalue of $A^\sigma$ associated with the
primitive idempotent $E^\sigma_i$. 
Similarly $\th^*_i$ is the eigenvalue of $A^{*\sigma}$ associated with the
primitive idempotent $E^{*\sigma}_i$. The result follows.

(iii):
The primitive idempotent $E_0$ has rank $1$ since $\Phi$ is sharp,
so $\tr(E_0)=1$. By this and Lemma \ref{lem:antitr} we have 
$\tr(E_0^\sigma)=1$. 
Now the primitive idempotent $E_0^\sigma$ has rank $1$ so $\Phi^\sigma$ is sharp.
By \eqref{eq:zetai1}, Lemma \ref{lem:antitr}, and (ii) above,
we find $\Phi$ and $\Phi^\sigma$
have the same split sequence.
\end{proof}

\medskip

Let $\tilde{V}$ denote the dual space of $V$; consisting of
all $\mathbb{K}$-linear transformations from $V$ to $\mathbb{K}$.
Define $\b{\;,\,} : V \times \tilde{V} \to \mathbb{K}$ by
$\b{v,f}=f(v)$ for $v \in V$, $f \in \tilde{V}$.
Then $\b{\;,\,}$ is a nondegenerate bilinear form.
We call this form
the {\em canonical bilinear form between $V$ and $\tilde{V}$}.
Let $\sigma : \text{End}(V) \to \text{End}(\tilde{V})$
denote the anti-isomorphism associated with $\b{\;,\,}$,
so that $(X^\sigma f)v=f(Xv)$ for all $v \in V$, $f \in \tilde{V}$,
$X \in \text{End}(V)$.
We call $\sigma$ the {\em canonical anti-isomorphism} from $\text{End}(V)$
to $\text{End}(\tilde{V})$.

\medskip

\begin{definition}  {\rm \cite{AC3}}     \label{def:dualTDS}  \samepage
Let $\Phi=(A;\{E_i\}_{i=0}^d;A^*;\{E^*_i\}_{i=0}^d)$ denote
a tridiagonal system on $V$.
Let $\tilde{V}$ denote the dual space of $V$ and let
$\sigma : \text{End}(V) \to \text{End}(\tilde{V})$ denote the
canonical anti-isomorphism.
By Proposition \ref{prop:antiTD}, $\Phi^\sigma$ is a tridiagonal system on 
$\tilde{V}$; we call this tridiagonal system the {\em dual} of $\Phi$.
\end{definition}

\begin{corollary}           \label{cor:dualPhi}    \samepage
With reference to Definition {\rm \ref{def}} assume $\Phi$ is sharp. 
Then the dual of $\Phi$ is sharp and has the same parameter array
as $\Phi$.
\end{corollary}

\begin{proof}
Immediate from Proposition \ref{prop:antiTD}
and Definition \ref{def:dualTDS}.
\end{proof}

\begin{proposition}    \label{prop:antibilin}  \samepage
With reference to Definition {\rm \ref{def}}
let $\b{\;,\,}$ denote a nondegenerate bilinear form on $V$ and let
$\dagger$ denote the associated anti-automorphism of $\text{\rm End}(V)$. 
Then the following {\rm (i)}, {\rm (ii)} are equivalent.
\begin{itemize}
\item[\rm (i)] 
The bilinear form $\b{\;,\,}$ satisfies \eqref{eq:compatible}.
\item[\rm (ii)]
$A^\dagger = A$ and ${A^*}^\dagger = A^*$. 
\end{itemize}
\end{proposition}

\begin{proof}
(i)$\Rightarrow$(ii):
Concerning $A$, we fix $v \in V$ and show $(A-A^\dagger)v=0$. 
For all $u \in V$ we have $\b{Au,v} = \b{u,A^\dagger v}$  
by construction and
$\b{Au,v} = \b{u,Av}$ by \eqref{eq:compatible}. 
Therefore $\b{u, (A-A^\dagger)v} = 0$ and this gives
$(A-A^\dagger)v = 0$ since $\b{\;,\,}$ is nondegenerate.
We have shown $A^\dagger = A$ and one similarly shows
${A^*}^\dagger = A^*$.

(ii)$\Rightarrow$(i):
For $u,v \in V$, by construction $\b{Xu,v}= \b{u, X^\dagger v}$
for all $X \in \text{End}(V)$. 
In this equation we set $X=A$ and use $A^\dagger = A$ to get
$\b{Au,v} = \b{u, Av}$.
Similarly we obtain $\b{A^*u,v} = \b{u, A^*v}$. 
Therefore $\b{\;,\,}$ satisfies \eqref{eq:compatible}.
\end{proof}

\section{Conjecture \ref{conj:main} and bilinear forms}

\indent
In this section we show how the truth of Conjecture \ref{conj:main}
implies the existence of a nonzero bilinear form 
that satisfies \eqref{eq:compatible}.
Actually, to obtain this bilinear form 
we do not need the full strength 
of Conjecture \ref{conj:main}; just the following weaker version.

\medskip

\begin{conjecture}           \label{conj}  \samepage
Two sharp tridiagonal systems over $\mathbb{K}$
are isomorphic if and only if they have the same parameter array.
\end{conjecture}

\begin{note}
In \cite{IT:aug} Conjecture \ref{conj} is proved assuming
the parameter $q$ from above Lemma \ref{lem:formula} is not a root of unity,
and that $\mathbb{K}$ is algebraically closed.
In \cite{IT:Krawt} Conjecture \ref{conj} is proved assuming the tridiagonal
pairs in question have Krawtchouk type, and that $\mathbb{K}$ is algebraically
closed.
\end{note}

\begin{corollary}         \label{cor:dualisom}  \samepage
Assume Conjecture {\rm \ref{conj}} is true.
Then every sharp tridiagonal system over $\mathbb{K}$ is isomorphic to its dual.
\end{corollary}

\begin{proof}
Follows from Corollary \ref{cor:dualPhi}.
\end{proof}

\begin{theorem}          \label{thm:bilin}  \samepage
Assume Conjecture {\rm \ref{conj}}  is true.
With reference to Definition {\rm \ref{def}} assume $\Phi$ is sharp.
Then there exists a nonzero bilinear form on $V$ that satisfies
\eqref{eq:compatible}.
This form is unique up to multiplication by a nonzero scalar in $\mathbb{K}$.
This form is nondegenerate and symmetric.
\end{theorem}

\begin{proof}
Let $\tilde{V}$ denote the dual space of $V$. 
Let $(\;,\,) : V \times \tilde{V} \to \mathbb{K}$ denote
the canonical bilinear form and let 
$\sigma : \text{End}(V) \to \text{End}(\tilde{V})$
denote the canonical anti-isomorphism.
By Lemma \ref{lem:isomalt} and Corollary \ref{cor:dualisom} there exists an
isomorphism of vector spaces $\gamma: V \to \tilde{V}$
such that $\gamma A= A^\sigma \gamma$
and  
 $\gamma A^*= {A^*}^\sigma \gamma$.
Define a map $\b{\;,\,} : V \times V \to \mathbb{K}$ by
\[
  \b{u,v} = (u, \gamma v)  \qquad\qquad           (u,v \in V)
\]
and observe that $\b{\;,\,}$ is a bilinear form on $V$.
The form $\b{\;,\,}$ is nondegenerate since $(\,,\,)$ is nondegenerate. 
We show $\b{\;,\,}$ satisfies \eqref{eq:compatible}. 
For $u,v \in V$,
\[
 \b{Au,v}=(Au,\gamma v) = (u, A^\sigma \gamma v)
   =(u,\gamma Av) = \b{u,Av}
\]
and similarly
$\b{A^*u,v} = \b{u,A^*v}$.
Therefore $\b{\;,\,}$ satisfies \eqref{eq:compatible}.
By Lemma \ref{lem:unique},
$\b{\;,\,}$ is unique up to multiplication by a nonzero scalar in $\mathbb{K}$.
By Lemma \ref{lem:symmetric}, $\b{\;,\,}$ is symmetric.
\end{proof}

\begin{theorem}           \label{thm:existanti}  \samepage
Assume Conjecture {\rm \ref{conj}} is true.
With reference to Definition {\rm \ref{def}} assume $\Phi$ is sharp.
Then there exists a unique anti-automorphism $\dagger$ of $\text{\rm End}(V)$
that fixes each of $A,A^*$.
Moreover $X^{\dagger\dagger}=X$ for all $X \in \text{\rm End}(V)$.
\end{theorem}

\begin{proof}
Let $\b{\;,\,}$ denote the bilinear form on $V$ from Theorem \ref{thm:bilin},
and let $\dagger$ denote the associated anti-automorphism
of $\text{End}(V)$. 
Then $A^\dagger = A$ and ${A^*}^\dagger = A^*$ 
by Proposition \ref{prop:antibilin}.
The anti-automorphism $\dagger$ is unique by Lemma \ref{lem:unique} and 
Proposition \ref{prop:antibilin}.
For $X \in \text{End}(V)$ we have $X^{\dagger\dagger} = X$ by 
Lemma \ref{lem:involutive} and
since $\b{\;,\,}$ is symmetric by Lemma \ref{lem:symmetric}.
\end{proof}

\section{Directions for further research}

\indent
In this section we give some suggestions for further research.
We start with a definition.

\medskip

\begin{definition}    \label{def:x}   \samepage
With reference to Definition \ref{def} let $\cal D$
(resp. ${\cal D}^*$) denote the $\mathbb{K}$-subalgebra of 
$\text{End}(V)$ generated by $A$ (resp. $A^*$).
Let $T$ denote the $\mathbb{K}$-subalgebra of $\text{End}(V)$
generated by $A$ and $A^*$.
\end{definition}

\medskip

We use the following convention.
For $\mathbb{K}$-subspaces $R,S$ of $\text{\rm End}(V)$ let
$RS$ denote the $\mathbb{K}$-subspace spanned by
$\{rs\,|\, r \in R, s \in S\}$.

\medskip

\begin{conjecture}   \label{conj:EDDE} \samepage
With reference to Definition {\rm \ref{def}} and Definition {\rm \ref{def:x}} 
the following {\rm (i)}--{\rm (iv)} hold.
\begin{itemize}
\item[\rm (i)]
The elements of $E^*_0 {\cal D} E^*_0$ mutually commute.
\item[\rm (ii)]
Each of the following holds:
\begin{align*}
E^*_0 {\cal D} {\cal D}^* E_0  
    &= E^*_0 {\cal D} E^*_0 E_0,   \\
E^*_0 {\cal D} {\cal D}^* {\cal D} E^*_0 
    &=  E^*_0 {\cal D} E^*_0 {\cal D} E^*_0,   \\
E^*_0 {\cal D} {\cal D}^* {\cal D} {\cal D}^* E_0 
    &=  E^*_0 {\cal D} E^*_0 {\cal D} E^*_0 E_0,    \\ 
E^*_0 {\cal D} {\cal D}^* {\cal D} {\cal D}^* {\cal D} E^*_0 
    &= E^*_0 {\cal D} E^*_0 {\cal D} E^*_0 {\cal D} E^*_0,    \\
E^*_0 {\cal D} {\cal D}^* {\cal D} {\cal D}^* {\cal D} {\cal D}^* E_0 
    &= E^*_0 {\cal D} E^*_0 {\cal D} E^*_0 {\cal D} E^*_0 E_0, \\
 \ldots & \ldots\ldots \\
 \ldots & \ldots\ldots
\end{align*}
\item[\rm (iii)]
The $\mathbb{K}$-algebra $E^*_0TE^*_0$ is generated by $E^*_0{\cal D}E^*_0$.
\item[\rm (iv)]
The $\mathbb{K}$-algebra $E^*_0TE^*_0$ is commutative.
\end{itemize}
\end{conjecture}

\begin{note}
In \cite{IT:aug}  Conjecture \ref{conj:EDDE} is proven 
assuming the parameter $q$ from above Lemma \ref{lem:formula} is not 
a root of unity,
and that $\mathbb{K}$ is algebraically closed.
\end{note}

\begin{conjecture}      \label{conj:field}   \samepage
With reference to Definition {\rm \ref{def}} and Definition {\rm \ref{def:x}} 
the following {\rm (i)}, {\rm (ii)} hold.
\begin{itemize}
\item[\rm (i)]
$E^*_0TE^*_0$ is a field with identity $E^*_0$.
\item[\rm (ii)]
Viewing $\mathbb{K}E^*_0$ as a field with identity $E^*_0$, the field
$E^*_0TE^*_0$ is an $r$-dimensional field extension of 
$\mathbb{K}E^*_0$, where $r = \dim E^*_0V$.
\end{itemize}
\end{conjecture}

\begin{note}   \samepage
Conjecture \ref{conj:field} implies Conjecture \ref{conj:closed}, since if
$\mathbb{K}$ is algebraically closed then the field 
$\mathbb{K}E^*_0$ has no finite-dimensional field
extensions other than itself.
\end{note}

\begin{problem}          \samepage
With reference to Definition \ref{def},
let $\{U_i\}_{i=0}^d$ denote the $\Phi$-split decomposition of $V$
from Section 5.
By \eqref{eq:up} and \eqref{eq:down},
for $0 \leq i \leq d/2$ the space $U_i$ is invariant under
\begin{equation}    \label{eq:invUi}
 (A^*-\th^*_{i+1}I)(A^*-\th^*_{i+2}I) \cdots (A^*- \th^*_{d-i}I)
  (A-\th_{d-i-1}I) \cdots (A-\th_{i+1}I)(A- \th_i I).
\end{equation}
The restriction of \eqref{eq:invUi} to $U_i$ is invertible 
by \cite[Lemma 6.5]{ITT}.
Assuming $\Phi$ is sharp,
find the eigenvalues for the action of \eqref{eq:invUi} on $U_i$
in terms of the parameter array of $\Phi$.
\end{problem}

\begin{problem}   \samepage
With reference to Definition \ref{def}, let 
$\{U_i\}_{i=0}^d$ denote the $\Phi$-split decomposition of $V$ from Section 5.
For $0 \leq i \leq d$ define a linear transformation
$F_i:V\to V$ such that $(F_i-I)U_i=0$ and 
$F_iU_j=0$ if $i\neq j$ $(0 \leq j\leq d)$.
In other words $F_i$ is the projection from $V$ onto $U_i$.
Assuming $A,A^*$ is sharp, find $F_i$ as a polynomial in $A, A^*$ and
express the coefficients of this polynomial in terms
of the parameter array.
\end{problem}

\begin{problem}     \samepage
With reference to Definition \ref{def} assume $\Phi$ is sharp.
For $0 \leq i \leq d$ find each of
\[
 \tr(E_iE^*_0),  \qquad  \tr(E_iE^*_d), \qquad
 \tr(E^*_iE_0),  \qquad  \tr(E^*_iE_d)
\]
in terms of the parameter array of $\Phi$.
\end{problem}

\bigskip

{\small

\bibliographystyle{plain}

\begin{thebibliography}{1}

\bibitem{AC}
H.~Alnajjar, B.~Curtin,
A family of tridiagonal pairs,
Linear Algebra Appl.\ 390 (2004) 369--384.

\bibitem{AC2}
H.~Alnajjar, B.~Curtin,
A family of tridiagonal pairs related to the quantum affine 
algebra $U\sb q(\widehat{\mathfrak{sl}}\sb 2)$,
Electron. J. Linear Algebra 13 (2005) 1--9. 

\bibitem{AC3}
H.~Alnajjar, B.~Curtin,
A bilinear form for tridiagonal pairs of $q$-Serre type,
submitted for publication.

\bibitem{BI}
E. Bannai, T. Ito,
Algebraic Combinatorics I: Association Schemems,
Benjamin/Cummings, London 1984.

\bibitem{Bas}
P.~Baseilhac,
A family of tridiagonal pairs and related symmetric functions,
J. Phys. A  39 (38) (2006) 11773--11791.

\bibitem{BT:Borel}
G.~Benkart, P.~Terwilliger,
Irreducible modules for the quantum affine algebra
$U_q(\widehat{\mathfrak{sl}}_2)$ and its Borel subalgebra,
J. Algebra 282 (2004) 172--194;
{\tt arXiv:math.QA/0311152}.

\bibitem{BT:loop}
G.~Benkart, P.~Terwilliger,
The universal central extension of the three-point
$\mathfrak{sl}_2$ loop algebra,
Proc.\ Amer.\ Math.\ Soc.\ 135 (6) (2007) 1659--1668;
{\tt arXiv:math.RA/0512422}.

\bibitem{Bow}
J.~Bowman,
Irreducible modules for the quantum affine algebra
$U_q(\mathfrak{g})$ and its Borel subalgebra $U_q(\mathfrak{g})^{\geq 0}$,
preprint;
{\tt arxiv:math.QA/0606627}.
       
\bibitem{Ca}
J.S.~Caughman,
The Terwilliger algebras of bipartite $P$- and $Q$-polynomial schemes,
Discrete Math. 196 (1999) 65--95.

\bibitem{CaMT}
J.S.~Caughman, M.S.~MacLean, P.~Terwilliger,
The Terwilliger algebra of an almost-bipartite $P$- and $Q$-polynomial 
association scheme,
Discrete Math. 292 (2005) 17--44.

\bibitem{CaW}
J.S.~Caughman, N.~Wolff,
The Terwilliger algebra of a distance-regular graph 
that supports a spin model,
J. Algebraic Combin. 21 (2005) 289--310.

\bibitem{Cur:mlt}
B.~Curtin,
Modular Leonard triples,
Linear Algebra Appl. 424 (2007) 510--539.

\bibitem{Cur:spinLP}
B.~Curtin,
Spin Leonard pairs, 
Ramanujan J.\ 13 (2007) 319--332.

\bibitem{Egge}
E.~Egge,
A generalization of the Terwilliger algebra,
J. Algebra 233 (2000) 213--252.

\bibitem{F:RL}
D.~Funk-Neubauer,
Raising/lowering maps and modules for the
quantum affine algebra $U_q(\widehat{\mathfrak{sl}}_2)$,
Comm.\ Algebra 35 (2007) 2140--2159;
{\tt arXiv:math.QA/0506629}.


\bibitem{H}
B.~Hartwig,
Three mutually adjacent Leonard pairs,
Linear Algebra Appl. 408 (2005) 19--39;
{\tt arXiv:math.AC/0508415}.

\bibitem{H:tetra}
B.~Hartwig,
The Tetrahedron algebra and its finite dimensional irreducible modules,
Linear Algebra Appl. 422 (2007) 219--235;
{\tt arXiv:math.RT/0606197}.

\bibitem{HT:tetra}
B.~Hartwig, P.~Terwilliger,
The tetrahedron algebra, the Onsager algebra, 
and the $\mathfrak{sl}_2$ loop algebra,
J. Algebra 308 (2007) 840--863;
{\tt arXiv:math-ph/0511004}.


\bibitem{ITT}
T.~Ito, K.~Tanabe, P.~Terwilliger,
Some algebra related to ${P}$- and ${Q}$-polynomial association schemes, 
Codes and Association Schemes (Piscataway NJ, 1999), 
American Mathematical Society, Providence, RI, 2001, pp. 167--192;
{\tt arXiv:math.CO/0406556}.

\bibitem{IT:shape}
T.~Ito, P.~Terwilliger,
The shape of a tridiagonal pair,
J. Pure Appl.\ Algebra 188 (2004) 145--160;
{\tt arXiv:math.QA/0304244}.

\bibitem{IT:uqsl2hat}
T.~Ito, P.~Terwilliger,
Tridiagonal pairs and the quantum affine algebra
$U_q(\widehat{\mathfrak{sl}}_2)$,
Ramanujan J. 13 (2007) 39--62;
{\tt arXiv:math.QA/0310042}.

\bibitem{IT:non-nilpotent}
T.~Ito, P.~Terwilliger,
Two non-nilpotent linear transformations that satisfy the cubic $q$-Serre relations,
J. Algebra Appl. 6 (2007) 477--503;
{\tt arXiv:math.QA/0508398}.

\bibitem{IT:tetra}
T.~Ito, P.~Terwilliger,
The $q$-tetrahedron algebra and its finite dimensional irreducible modules,
Comm.\ Algebra, in press;
{\tt arXiv:math.QA/0602199}.

\bibitem{IT:inverting}
T.~Ito, P.~Terwilliger,
$q$-inverting pairs of linear transformations and the $q$-tetrahedron algebra,
Linear Algebra Appl. 427 (2007) 218--233;
{\tt arXiv:math.RT/0606237}.

\bibitem{IT:drg}
T.~Ito, P.~Terwilliger,
Distance-regular graphs and the $q$-tetrahedron algebra,
European J. Combin.,
submitted for publication;
{\tt arxiv:math.CO/0608694}.

\bibitem{IT:loop}
T.~Ito, P.~Terwilliger,
Finite-dimensional irreducible
modules for the three-point $\mathfrak{sl}_2$ loop algebra,
Comm.\ Algebra, submitted for publication;
{\tt arXiv:0707.2313}. 

\bibitem{IT:Krawt}
T.~Ito, P.~Terwilliger,
Tridiagonal pairs of Krawtchouk type,
Linear Algebra Appl., in press;
{\tt arXiv:0706.1065}.

\bibitem{IT:aug}
T.~Ito, P.~Terwilliger,
The augmented tridiagonal algebra,
preprint.

\bibitem{ITW:equitable}
T.~Ito, P.~Terwilliger, C.~Weng,
The quantum algebra $U_q(\mathfrak{sl}_2)$ and its equitable presentation,
J. Algebra 298 (2006) 284--301;
{\tt arXiv:math.QA/0507477}.

\bibitem{Koe}
R.~Koekoek, R.~F.~Swarttouw,
The Askey-scheme of hypergeometric orthogonal polynomials and its q-analogue,
report 98-17,
Delft University of Technology,
The Netherlands, 1998;
available at
{\tt http://aw.twi.tudelft.nl/~koekoek/askey.html}.

\bibitem{Leo}
D.~Leonard,
Orthogonal polynomials, duality, and association schemes,
SIAM J. Math. Anal. 13 (1982) 656--663.

\bibitem{Mik}
S.~Miklavic,
On bipartite Q-polynomial distance-regular graphs,
Europian J. Combin. 28 (2007) 94--110.

\bibitem{Mik2}
S.~Miklavic,
$Q$-polynomial distance-regular graphs with $a_1=0$ and $a_2 \not=0$,
submitted for publication.

\bibitem{M:LT}
S.~Miklavic,
Leonard triples and hypercubes,
submitted for publication.

\bibitem{N:aw}
K.~Nomura,
\newblock
Tridiagonal pairs and the {A}skey-{W}ilson relations,
Linear Algebra Appl.\ 397 (2005) 99--106.

\bibitem{N:refine}
K.~Nomura,
A refinement of the split decomposition of a tridiagonal pair,
Linear Algebra Appl.\ 403 (2005) 1--23.

\bibitem{N:height1}
K.~Nomura,
Tridiagonal pairs of height one,
Linear Algebra Appl.\ 403 (2005) 118--142.

\bibitem{NT:balanced}
K.~Nomura, P.~Terwilliger,
Balanced Leonard pairs,
Linear Algebra Appl. 420 (2007) 51--69.
{\tt arXiv:math.RA/0506219}. 

\bibitem{NT:formula}
K.~Nomura, P.~Terwilliger,
Some trace formulae involving the split sequences of a Leonard pair,
Linear Algebra Appl.\ 413 (2006) 189--201;
{\tt arXiv:math.RA/0508407}.

\bibitem{NT:det}
K.~Nomura, P.~Terwilliger,
The determinant of $AA^*-A^*A$ for a Leonard pair $A,A^*$,
Linear Algebra Appl.\ 416 (2006) 880--889;
{\tt arXiv:math.RA/0511641}.

\bibitem{NT:mu}
K.~Nomura, P.~Terwilliger,
Matrix units associated with the split basis of a Leonard pair,
Linear Algebra Appl.\ 418 (2006) 775--787;
{\tt arXiv:math.RA/0602416}.

\bibitem{NT:span}
K.~Nomura, P.~Terwilliger,
Linear transformations that are tridiagonal with respect to
both eigenbases of a Leonard pair,
Linear Algebra Appl. 420 (2007) 198--207;
{\tt arXiv:math.RA/0605316}.

\bibitem{NT:switch}
K.~Nomura, P.~Terwilliger,
The switching element for a Leonard pair,
Linear Algebra Appl., in press;
{\tt arXiv:math.RA/0608623}.

\bibitem{NT:affine}
K.~Nomura, P.~Terwilliger,
Affine transformations of a Leonard pair,
submitted for publication;
{\tt arXiv:math.RA/0611783}.

\bibitem{NT:tde}
K.~Nomura, P.~Terwilliger,
The split decomposition of a tridiagonal pair,
Linear Algebra Appl. 424 (2007) 339--345;
{\tt arXiv:math.RA/0612460}.

\bibitem{NT:maps}
K.~Nomura, P.~Terwilliger,
Transition maps between the 24 bases for a Leonard pair,
submitted for publication;
{\tt arXiv:0705.3918}.

\bibitem{P}
A.~A.~Pascasio,
On the multiplicities of the primitive idempotents of a
{$Q$}-polynomial distance-regular graph,
European J. Combin.\ 23 (2002) 1073--1078.

\bibitem{R:multi}
H.~Rosengren,
Multivariable orthogonal polynomials and coupling coefficients for 
discrete series representations,
SIAM J. Math. Anal.\ 30 (1999) 233--272.

\bibitem{R:6j}
H.~Rosengren,
An elementary approach to the $6j$-symbols
(classical, quantum, rational, trigonometric, and elliptic),
Ramanujan J. 13 (2007) 131--166;
{\tt arXiv:math.CA/0312310}.

\bibitem{Rot}
J.~J.~Rotman,
Advanced modern algebra,
Prentice Hall,
Saddle River NJ 2002.

\bibitem{Suz}
H.~Suzuki,
On strongly closed subgraphs with diameter two and the Q-polynomial property,
Europian J. Combin. 28 (2007) 167--185.

\bibitem{T:sub1}
P.~Terwilliger,
The subconstituent algebra of an association scheme I,
J. Algebraic Combin.\ 1 (1992) 363--388.

\bibitem{T:sub3}
P.~Terwilliger,
The subconstituent algebra of an association scheme III,
J. Algebraic Combin.\ 2 (1993) 177--210.

\bibitem{T:Leonard}
P.~Terwilliger,
Two linear transformations each tridiagonal with respect to an 
eigenbasis of the other,
Linear Algebra Appl.\ 330 (2001) 149--203;
{\tt arXiv:math.RA/0406555}.

\bibitem{T:qSerre}
P.~Terwilliger,
Two relations that generalize the $q$-Serre relations and the
Dolan-Grady relations,
Physics and Combinatorics 1999 (Nagoya), World Scientific Publishing,
River Edge, NJ, 2001, pp.\ 377--398;
{\tt arXiv:math.QA/0307016}.

\bibitem{T:24points}
P.~Terwilliger,
Leonard pairs from 24 points of view,
Rocky Mountain J. Math.\ 32 (2) (2002) 827--888;
{\tt arXiv:math.RA/0406577}.

\bibitem{T:canform}
P.~Terwilliger,
Two linear transformations each tridiagonal with respect to an
eigenbasis of the other; the $TD$-$D$ and the $LB$-$UB$ canonical form,
J. Algebra 291 (2005) 1--45;
{\tt arXiv:math.RA/0304077}.

\bibitem{T:intro}
P.~Terwilliger,
Introduction to Leonard pairs,
J.\ Comput.\ Appl.\ Math.\ 153 (2) (2003) 463--475.

\bibitem{T:intro2}
P.~Terwilliger,
Introduction to {L}eonard pairs and {L}eonard systems,
S\=uri\-kaiseki\-kenky\=usho K\=oky\=uroku 1109 (1999) 67--79, 
Algebraic Combinatorics (Kyoto, 1999).

\bibitem{T:split}
P.~Terwilliger,
Two linear transformations each tridiagonal with respect to an
eigenbasis of the other; comments on the split decomposition,
J. Comput.\ Appl.\ Math.\ 178 (2005) 437--452;
{\tt arXiv:math.RA/0306290}.

\bibitem{T:array}
P.~Terwilliger,
Two linear transformations each tridiagonal with respect to an
eigenbasis of the other; comments on the parameter array,
Des.\ Codes Cryptogr.\ 34 (2005) 307--332;
{\tt arXiv:math.RA/0306291}.

\bibitem{T:qRacah}
P.~Terwilliger,
Leonard pairs and the $q$-Racah polynomials,
Linear Algebra Appl.\ 387 (2004) 235--276;
{\tt arXiv:math.QA/0306301}.

\bibitem{T:survey}
P.~Terwilliger,
An algebraic approach to the Askey scheme of orthogonal polynomials,
in: Orthogonal Polynomials and Special Functions: computation and applications,
Lecture Notes in Mathematics, vol.\ 1883, Springer, 2006, pp.\ 255--330;
{\tt arXiv:math.QA/0408390}.

\bibitem{T:Kac-Moody}
P.~Terwilliger,
The equitable presentation for the quantum group $U_q(\mathfrak{g})$ 
associated with a symmetrizable Kac-Moody algebra $\mathfrak{g}$,
J. Algebra 298 (2006) 302--319; 
{\tt arXiv:math.QA/0507478}.

\bibitem{TV}
P.~Terwilliger, R.~Vidunas,
Leonard pairs and the Askey-Wilson relations,
J. Algebra Appl.\ 3 (2004) 411--426;
{\tt arXiv:math.QA/0305356}.

\bibitem{Vidar}
M.~Vidar,
Tridiagonal pairs of shape $(1,2,1)$, preprint.

\bibitem{V}
R.~Vidunas,
Normalized Leonard pairs and Askey-Wilson relations,
Technical Report MHF 2005-16, Kyushu University 2005;
{\tt arXiv:math.RA/0505041}.

\bibitem{V:AW}
R.~Vidunas,
Askey-Wilson relations and Leonard pairs,
Technical Report MHF 2005-17, Kyushu University 2005;
{\tt arXiv:math.QA/0511509}.

\bibitem{Z}
A.S.~Zhedanov,
``Hidden symmetry'' of Askey-Wilson polynomials,
Teoret.\ Mat Fiz.\ 89 (1991) 190--204. 
\end{thebibliography}

}

\bigskip\bigskip\noindent
Kazumasa Nomura\\
College of Liberal Arts and Sciences\\
Tokyo Medical and Dental University\\
Kohnodai, Ichikawa, 272-0827 Japan\\
email: knomura@pop11.odn.ne.jp

\bigskip\noindent
Paul Terwilliger\\
Department of Mathematics\\
University of Wisconsin\\
480 Lincoln Drive\\ 
Madison, Wisconsin, 53706 USA\\
email: terwilli@math.wisc.edu

\bigskip\noindent
{\bf Keywords.}
Leonard pair, tridiagonal pair, $q$-Racah polynomial, orthogonal polynomial.

\noindent
{\bf 2000 Mathematics Subject Classification}.
05E35, 05E30, 33C45, 33D45.

\end{document}